\documentclass[11pt]{article}
\usepackage{epsfig}

\def\be{\begin{equation}}
\def\ee{\end{equation}}
\def\bea{\begin{eqnarray}}
\def\eea{\end{eqnarray}}
\def\bes{\begin{eqnarray*}}
\def\ees{\end{eqnarray*}}

\def\nn{\nonumber}
\def\<{\langle}
\def\>{\rangle}
\def\lb{\label}

\def\R{{\bf R}}
\def\C{{\bf C}}
\def\Z{{\bf Z}}
\def\N{{\bf N}}
\def\U{{\bf U}}
\def\Q{{\bf Q}}

\def\bb{{\beta}}
\def\ga{{\gamma}}

\def\ka{{\kappa}}
\def\th{{\theta}}

\def\lm{{\lambda}}
\def\Lm{{\Lambda}}

\def\rank{{\rm rank}}
\def\Sp{{\rm Sp}}
\def\mod{{\rm mod}}
\def\grad{{\rm grad}}
\def\ol#1{\overline{#1}}  

\def\hb{\vrule height0.18cm width0.14cm $\,$}

\def\ol#1{\overline{#1}}  

\title{ Closed geodesics on positively curved Finsler spheres}
\author{Wei Wang\thanks{Partially supported by LMAM in Peking University
in China and China Postdoctoral Science Foundation No.20070420264.
E-mail: alexanderweiwang@yahoo.com.cn, wangwei@math.pku.edu.cn  }\\
School of Mathematical Science \\ Peking University, Beijing 100871 \\
PEOPLES REPUBLIC OF CHINA \\ }
\date{}

\begin{document}

\maketitle

\begin{abstract}
{\it In this paper, we prove that for every Finsler  $n$-sphere
$(S^n,\,F)$ for $n\ge 3$ with reversibility $\lambda$ and flag
curvature $K$ satisfying
$\left(\frac{\lambda}{\lambda+1}\right)^2<K\le 1$, either there
exist infinitely many prime closed geodesics or there exists one
elliptic closed geodesic whose linearized Poincar\'e map has at
least one eigenvalue which is of the form $\exp(\pi i \mu)$ with an
irrational $\mu$. Furthermore, there always exist three prime closed
geodesics on any $(S^3,\,F)$ satisfying the above pinching
condition.  }
\end{abstract}



{\bf Running head}: Closed geodesics on Finsler spheres

\renewcommand{\theequation}{\thesection.\arabic{equation}}
\renewcommand{\thefigure}{\thesection.\arabic{figure}}

\setcounter{equation}{0}
\section{Introduction and main results}

This paper is devoted to a study on closed geodesics on Finsler
$n$-spheres. For the definition of closed geodesics on a Finsler
manifold, we refer readers to \cite{BCS1} and \cite{She1}. As usual,
on any Finsler n-sphere $S^n=(S^n,F)$ a closed geodesic
$c:S^1=\R/\Z\to S^n$ is {\it prime}, if it is not a multiple
covering (i.e., iteration) of any other closed geodesics. Here the
$m$-th iteration $c^m$ of $c$ is defined by $c^m(t)=c(mt)$. The
inverse curve $c^{-1}$ of $c$ is defined by $c^{-1}(t)=c(1-t)$ for
$t\in \R$. Note that on a non-symmetric Finsler sphere, $c^{-1}$
is not a geodesic in general. We call two prime closed geodesics $c$ and $d$
{\it distinct}, if there is no $\th\in S^1$ such that
$c(t)=d(t+\th)$ for all $t\in\R$. We shall omit the word {\it
distinct} for short when we talk about more than one prime closed
geodesics. On a symmetric Finsler (or Riemannian) $n$-sphere, two closed geodesics
$c$ and $d$ are called { \it geometrically distinct}, if $
c(S^1)\neq d(S^1)$, i.e., their image sets in $S^n$ are distinct.

For a closed geodesic $c$ on $(S^n,\,F)$, denote by $P_c$
the linearized Poincar\'{e} map of $c$. Then $P_c\in \Sp(2n-2)$ is symplectic.
For any $M\in \Sp(2k)$, we define the {\it elliptic height } $e(M)$
of $M$ to be the total algebraic multiplicity of all eigenvalues of
$M$ on the unit circle $\U=\{z\in\C|\; |z|=1\}$ in the complex plane
$\C$. Since $M$ is symplectic, $e(M)$ is even and $0\le e(M)\le 2k$.
A closed geodesic $c$ is called {\it elliptic} if $e(P_c)=2(n-1)$, i.e., all the
eigenvalues of $P_c$ locate on $\U$.  Following H-B. Rademacher in
\cite{Rad3}, the reversibility $\lambda=\lambda(M,\,F)$ of a compact
Finsler manifold $(M,\,F)$ is defied to be
$$\lambda:=\max\{F(-X)\,|\,X\in TM, \,F(X)=1\}\ge 1.$$

It was quite surprising when Katok \cite{Kat1} in 1973 found some non-symmetric
Finsler metrics on $S^n$ with only finitely many prime closed geodesics and all closed
geodesics are non-degenerate and elliptic. The smallest number of closed geodesics
that one obtains in these examples is $2n$ on $S^{2n}$ and $S^{2n-1}$ (cf. \cite{Zil1}).
Then it is an open question whether there are always at least $n$ prime closed
geodesics on any Finsler $n$-sphere (cf. p.156 of \cite{Zil1}).

The following are the main results in this paper:

{\bf Theorem 1.1.} {\it For every Finsler  $n$-sphere $(S^n,\,F)$
for $n\ge 3$ with reversibility $\lambda$ and flag curvature $K$
satisfying $\left(\frac{\lambda}{\lambda+1}\right)^2<K\le 1$, either
there exist infinitely many prime closed geodesics or there exists
one elliptic closed geodesic whose linearized Poincar\'e map has at
least one eigenvalue which is of the form $\exp(\pi i \mu)$ with an
irrational $\mu$.  The number $\mu$ could be called a Floquet
exponent.   }

{\bf Theorem 1.2.} {\it  Under the assumption of Theorem 1.1 and
suppose the number of prime closed geodesics is finite. Denote by
$\{P_{c_j}\}_{1\le j\le p}$ the linearized Poincar\'e maps of them.
Then there exist at least two pairs of eigenvalues of those
$P_{c_j}$s which are of the form $\exp(\pi i \mu)$ with $\mu$
being irrational. (the
two pairs of eigenvalues maybe belong to one $P_{c_j}$). }

{\bf Remark 1.3.} Note that on the standard Riemannian $n$-sphere of constant
curvature $1$, all geodesics are closed and their linearized
Poincar\'e maps are $I_{2n-2}$, i.e., the identity matrix in $\R^{2n-2}$
and then there exists no eigenvalue of them which is an irrational multiple
of $\pi$. Hence our above theorems describe a character of a Finsler sphere
which carries finitely many prime closed geodesics.

Note also that our definition of ellipticity is different
from that in \cite{BTZ1}, in which they call a closed geodesic $c$
is of {\it elliptic-parabolic type} if the linearized Poincar\'e map $P_c$ of
$c$ splits into two-dimensional rotations and a part whose eigenvalues
are $\pm 1$.

In \cite{BTZ1}, W. Ballmann, G. Thorbergsson and W. Ziller proved that
if a Riemannian manifold  $M=S^n$ satisfies $\frac{9}{16}\le K\le1$, there
exists a prime closed geodesic of elliptic-parabolic type on $M$.
In \cite{Rad4}, H-B. Rademacher proved that if a Finsler manifold
$M=S^n$ satisfies $\frac{9}{4}\frac{\lambda^2}{(\lambda+1)^2}< K\le1$ with $\lambda<2$,
there exists a prime closed geodesic of elliptic-parabolic type on $M$.
Comparing with their results, we can prove the following

{\bf Theorem 1.4.} {\it Under the assumption of Theorem 1.1
and suppose the number of prime closed geodesics is finite.
Then there  exists a closed geodesic $c$ (not necessarily to be prime)
of  elliptic-parabolic type in the sense of \cite{BTZ1}. }

Note that by Theorem 7 in \cite{Rad4}, the existence of at least two prime
closed geodesics on any Finsler  $3$-sphere $(S^3,\,F)$ satisfying
the assumption of Theorem 1.1. Now we can prove the following

{\bf Theorem 1.5.} {\it For every Finsler  $3$-sphere
$(S^3,\,F)$  with reversibility $\lambda$ and flag
curvature $K$ satisfying
$\left(\frac{\lambda}{\lambda+1}\right)^2< K\le 1$,  there
exist  at least three prime closed geodesics.    }

{\bf Remark 1.6.} Note that our Theorem 1.5 gives a partial result
to the problem (4) in p.156 in \cite{Zil1} for the $S^3$ case.
As mentioned above, there are examples with exactly four
prime closed geodesics on the $3$-sphere.
In \cite{BTZ2}, W. Ballmann, G. Thorbergsson and
W. Ziller proved that for a Riemannian metric on $S^n$ with sectional curvature $1/4\le K\le 1$ there exist
$g(n)$ geometrically distinct closed geodesics, and $g(3)=4$.
In \cite{LoW1}, Y. Long and the author proved that for every Riemannian
$3$-sphere $(M,g)$ with injectivity radius ${\rm inj}(M)\ge \pi$ and the sectional
curvature $K$ satisfying $\frac{1}{16}<K\le 1$ there exist at least two
geometrically distinct closed geodesics.
However, the proofs of these results relies essentially on the
symmetry of the metric, hence it can not be used here.
The method in \cite{LoW1} works only in the study of the
existence of two closed geodesics, while that in the present paper
works for more that two closed geodesics.

Our proof of these theorems contains mainly three ingredients: the
common index jump theorem of Y. Long, Morse theory, and
an existence theorem of N. Hingston. Fix a Finsler metric $F$
on $S^n$. Let $\Lm=\Lm S^n$ be the free loop space of $S^n$, which
is a Hilbert manifold. For definition and basic properties of
$\Lm$, we refer readers to \cite{Kli2} and \cite{Kli3}. Let
${E}(c)=\frac{1}{2}\int_0^1F(\dot{c}(t))^2dt$ be the energy
functional on $\Lm$. In this paper for $\ka\in \R$ we denote by
\be \Lm^{\ka}=\{d\in \Lm\,|\,E(d)\le \ka\},    \lb{1.1}\ee
and consider the quotient space $\Lm/S^1$. Since the
energy functional ${E}$ is $S^1$-invariant, the negative
gradient flow of $E$ induce a flow on $\Lm/S^1$, so we can apply
Morse theory on $\Lm/S^1$. By a result of H-B. Rademacher in
\cite{Rad1} of 1989, we get the Morse series of the space pair
$(\Lm/S^1, \Lm^0/S^1)$ with rational coefficients. The reason we
use $(\Lm/S^1, \Lm^0/S^1)$ instead of $(\Lm, \Lm^0)$ is that
the Morse series of the first is lacunary.

Sections 2 to 4 are preliminary materials for our proof.
In Section 2, basic properties of critical modules of closed
geodesics are introduced, whose proofs can be found in \cite{Rad2}
and \cite{BaL1}. In Section 3, Morse inequalities on the quotient space
$(\Lm/S^1,\Lm^0/S^1)$ are given, whose proof can be found in \cite{Rad1}
and \cite{Rad2}. In Section 4, by results in \cite{Lon1} of 2000
of Y. Long, we give the classification of closed geodesics on $S^n$.
N. Hingston's  Theorem in \cite{Hin1} is also listed in Section 4.

In Section 5, we establish a mean index equality when there exist only
finitely many prime closed geodesics on $(S^n,F)$. An abstract version of
such an equality was established by H-B. Rademacher in \cite{Rad2}.  Since
we need this equality with exact coefficients. we give a complete proof for
it in Section 5.

Based on these preparations, our Theorems 1.1-1.4 are proved in Section 6.
As an application of these theorems, Theorem 1.5 is proved in Section 7.

In this paper, let $\N$, $\N_0$, $\Z$, $\Q$, $\R$, and $\C$ denote
the sets of natural integers, non-negative integers, integers,
rational numbers, real numbers, and complex numbers respectively.
Let $(p,q)$ denotes the greatest common devisor of $p$ and $q\in\N$.
We use only singular homology modules with $\Q$-coefficients.
For terminologies in algebraic topology we refer to \cite{GrH1}.
For $k\in\N$, we denote by $\Q^k$ the direct sum $\Q\oplus\cdots\oplus\Q$ of
$k$ copies of $\Q$ and $\Q^0=0$. For an $S^1$-space $X$,
we denote by $\overline{X}$ the quotient space $X/S^1$.
We define the functions
\be \left\{\matrix{[a]=\max\{k\in\Z\,|\,k\le a\}, &
           \mathcal{E}(a)=\min\{k\in\Z\,|\,k\ge a\} , \cr
    \varphi(a)=\mathcal{E}(a)-[a],   \cr}\right. \lb{1.2}\ee
Especially, $\varphi(a)=0$ if $ a\in\Z\,$, and $\varphi(a)=1$ if $
a\notin\Z\,$.

\setcounter{equation}{0}
\section{Critical modules of iterations of closed geodesics}

In this section, we will study critical modules of closed geodesic,
all the details can be found in \cite{Rad2} or \cite{BaL1}.

On a compact Finsler manifold $(M,F)$, we choose an auxiliary Riemannian
metric. This endows the space $\Lambda=\Lambda M$ of $H^1$-maps
$\gamma:S^1\rightarrow M$ with a natural structure of Riemannian
Hilbert manifolds on which the group $S^1=\R/\Z$ acts continuously
by isometries, cf. \cite{Kli2}, Chapters 1 and 2. This action is
defined by translating the parameter, i.e.
$$ (s\cdot\gamma)(t)=\gamma(t+s) \qquad $$
for all $\gamma\in\Lm$ and $s,t\in S^1$.
The Finsler metric $F$ defines an energy functional $E$ and a length
functional $L$ on $\Lambda$ by
\be E(\gamma)=\frac{1}{2}\int_{S^1}F(\dot{\gamma}(t))^2dt,
 \quad L(\gamma) = \int_{S^1}F(\dot{\gamma}(t))dt.  \lb{2.1}\ee
Both functionals are invariant under the $S^1$-action. The critical points
of $E$ of positive energies are precisely the closed geodesics $c:S^1\to M$
of the Finsler structure. If $c\in\Lambda$ is a closed geodesic then $c$ is
a regular curve, i.e. $\dot{c}(t)\not= 0$ for all $t\in S^1$, and this implies
that the second differential $E''(c)$ of $E$ at $c$ exists.
As usual we define the index $i(c)$ of $c$ as the maximal dimension of
subspaces of $T_c \Lambda$ on which $E^{\prime\prime}(c)$ is negative definite, and the
nullity $\nu(c)$ of $c$ so that $\nu(c)+1$ is the dimension of the null
space of $E^{\prime\prime}(c)$.

For $m\in\N$ we denote the $m$-fold iteration map
$\phi^m:\Lambda\rightarrow\Lambda$ by \be \phi^m(\ga)(t)=\ga(mt)
\qquad \forall\,\ga\in\Lm, t\in S^1. \lb{2.2}\ee We also use the
notation $\phi^m(\gamma)=\gamma^m$. For a closed geodesic $c$, the
mean index is defined to be: \be
\hat{i}(c)=\lim_{m\rightarrow\infty}\frac{i(c^m)}{m}. \lb{2.3}\ee

If $\gamma\in\Lambda$ is not constant then the multiplicity
$m(\gamma)$ of $\gamma$ is the order of the isotropy group $\{s\in
S^1\mid s\cdot\gamma=\gamma\}$. If $m(\gamma)=1$ then $\gamma$ is
called {\it prime}. Hence $m(\gamma)=m$ if and only if there exists a
prime curve $\tilde{\gamma}\in\Lambda$ such that
$\gamma=\tilde{\gamma}^m$.

For a closed geodesic $c$ we set
$$ \Lm(c)=\{\ga\in\Lm\mid E(\ga)<E(c)\}. $$
If $A\subseteq\Lm$ is invariant under some subgroup $\Gamma$ of $S^1$,
we denote by $A/\Gamma$ the
quotient space of $A$ with respect to the action of $\Gamma$.

Using singular homology with rational coefficients we will
consider the following critical $\Q$-module of a closed geodesic
$c\in\Lambda$:
\be \overline{C}_*(E,c)
   = H_*\left((\Lm(c)\cup S^1\cdot c)/S^1,\Lm(c)/S^1\right). \lb{2.4}\ee

In order to relate the critical modules to the index and nullity
of $c$ we use the results of D. Gromoll and W. Meyer
from \cite{GrM1}, \cite{GrM2}. Following \cite{Rad2}, Section
6.2, we introduce finite-dimensional
approximations to $\Lambda$. We choose an arbitrary energy value
$a>0$ and $k\in\N$ such that every $F$-geodesic of length
$<\sqrt{2a/k}$ is minimal. Then
$$ \Lm(k,a)=\left\{\ga\in\Lm \mid E(\ga)<a \mbox{ and }
    \ga|_{[i/k,(i+1)/k]}\mbox{ is an $F$-geodesic for }i=0,\ldots,k-1\right\} $$
is a $(k\cdot\dim M)$-dimensional submanifold of $\Lambda$
consisting of closed geodesic polygons with $k$ vertices. The set
$\Lambda(k,a)$ is invariant under the subgroup $\Z_k$ of $S^1$.
Closed geodesics in $\Lambda^{a-}=\{\gamma\in\Lambda\mid
E(\gamma)<a\}$ are precisely the critical points of
$E|_{\Lm(k,a)}$, and for every closed geodesic $c\in\Lm(k,a)$ the
index of $(E|_{\Lm(k,a)})''(c)$ equals $i(c)$ and the null space
of $(E|_{\Lm(k,a)})''(c)$ coincides with the nullspace of
$E''(c)$, cf. \cite{Rad2}, p.51.

We call a closed geodesic satisfying the isolation condition, if
the following holds:

{\bf (Iso)  For all $m\in\N$ the orbit $S^1\cdot c^m$ is an
isolated critical orbit of $E$. }

Note that if the number of prime closed geodesics on a Finsler manifold
is finite, then all the closed geodesics satisfy (Iso).

Now we can apply the results by D. Gromoll and W. Meyer
\cite{GrM1} to a given closed geodesic $c$ satisfying (Iso). If
$m=m(c)$ is the multiplicity of $c$, we choose a finite-dimensinal
approximation $\Lm(k,a)\subseteq\Lm$ containing $c$ such that $m$
divides $k$. Then the isotropy subgroup $\Z_m\subseteq S^1$ of $c$
acts on $\Lm(k,a)$ by isometries.  Let $D$ be a $\Z_m$-invariant
local hypersurface transverse to $S^1\cdot c$ in $c\in D$.
According to \cite{GrM1}, Lemma 1, for every such $D$ we can find a product
neighborhood $B_+\times
B_-\times B_0$ of $0\in\R^{\dim\Lm(k,a)-1}$ such that $\dim
B_-=i(c), \dim B_0=\nu(c)$, and a diffeomorphism
$$  \psi:B=B_+\times B_-\times B_0\rightarrow\psi(B_+\times B_-\times B_0)\subseteq D  $$
from $B$ onto an open subset $\psi(B)\subseteq D$ such that
$\psi(0)=c$ and $\psi$ is $\Z_m$-invariant, and there exists a
smooth function $f:B_0\rightarrow\R$ satisfying \be f^\prime(0)=0
\quad \mbox{ and }\quad f^{\prime\prime}(0)=0  \lb{2.5}\ee and \be
E\circ\psi(x_+,x_-,x_0)=|x_+|^2 - |x_-|^2 + f(x_0),  \lb{2.6}\ee
for $(x_+, x_-, x_0) \in B_+\times B_-\times B_0$. As usual, we call
$$ N=\{\psi(0,0,x_0)|x_0\in B_0\} $$
a  local characteristic manifold at $c$ and
$$  U=\{\psi(0,x_-,0)|x_-\in B_-\}  $$
a  local negative disk  at $c$, $N$ and $U$ are
$\Z_m$-invariant. It follows from (\ref{2.6}) that $c$ is an
isolated critical point of $E|_N$. We set $N^-=N\cap\Lm(c)$,
$U^-=U\cap\Lm(c)=U\setminus\{c\}$ and $D^-=D\cap\Lm(c)$. Using
(\ref{2.6}), the fact that $c$ is an isolated critical point
of $E|_N$, and the K\"unneth formula, one concludes
\be H_*(D^-\cup\{c\},D^-)=
      H_*(U^-\cup\{c\},U^-) \otimes H_*(N^-\cup\{c\},N^-), \lb{2.7}\ee
where \be H_q(U^-\cup\{c\}, U^-) = H_q(U, U\setminus\{c\})
    = \left\{\matrix{\Q, & {\rm if\;}q=i(c), \cr
                      0, & {\rm otherwise}, \cr}\right.  \lb{2.8}\ee
cf. \cite{Rad2}, Lemma 6.4 and its proof. As studied in p.59 of \cite{Rad2},
for all $m\in\N$, let respectively
\be H_{\ast}(X,A)^{\pm\Z_m}
   = \{[\xi]\in H_{\ast}(X,A)\,|\,T_{\ast}[\xi]=\pm [\xi]\}, \lb{2.9}\ee
where $T$ is a generator of the $\Z_m$-action.

The following Propositions were proved in \cite{Rad2} and \cite{BaL1}.

{\bf Proposition 2.1.} (cf. Satz 6.11 of \cite{Rad2} or Proposition
3.12 of \cite{BaL1}) {\it Let $c$ be a prime closed geodesic on a Finsler
manifold $(M,F)$ satisfying (Iso). Then we have
\bea \overline{C}_q( E,c^m)
&&\equiv H_q\left((\Lm(c^m)\cup S^1\cdot c^m)/S^1, \Lm(c^m)/S^1\right)\nn\\
&&= \left(H_{i(c^m)}(U_{c^m}^-\cup\{c^m\},U_{c^m}^-)
    \otimes H_{q-i(c^m)}(N_{c^m}^-\cup\{c^m\},N_{c^m}^-)\right)^{+\Z_m} \nn
\eea

(i) When $\nu(c^m)=0$, there holds
\bea \overline{C}_q( E,c^m) = \left\{\matrix{
     \Q, &\quad {\it if}\;\; i(c^m)-i(c)\in 2\Z,\;{\it and}\;
                   q=i(c^m)\;  \cr
     0, &\quad {\it otherwise}. \cr}\right.  \lb{2.10}\eea

(ii) When $\nu(c^m)>0$, there holds
\bea \overline{C}_q( E,c^m) =
H_{q-i(c^m)}(N_{c^m}^-\cup\{c^m\}, N_{c^m}^-)^{\bb(c^m)\Z_m},\lb{2.11}\eea
where $\bb(c^m)=1$ (or $=-1$), when $i(c^m)-i(c)$ is even (or odd).}

In order to study the degenerate part of the local critical
modules, we need the following result which follows from Satz 6.6
of \cite{Rad2} directly.

{\bf Lemma 2.2.} {\it Let $c$ be a prime closed geodesic on a
Finsler manifold $(M, F)$ satisfying (Iso). Then there holds }
$$ H_q(N_{c^m}^-\cup\{c^m\}, N_{c^m}^-)^{+\Z_m}
 = H_q(N_{c^m}^-\cup\{c^m\}/\Z_m, N_{c^m}^-/\Z_m),
  \qquad \forall\,q\in\Z,\,m\in\N. $$

We introduce the following

{\bf Definition 2.3.} {\it Suppose $c$ is a closed geodesic of
multiplicity $m(c)=m$ satisfying (Iso). If $N$ is a local
characteristic manifold at $c$, $N^- = N\cap\Lm(c)$ and $j\in\Z$,
we define
\bea
k_j(c) &\equiv& \dim\, H_j( N^-\cup\{c\},N^-),    \lb{2.12}\\
k_j^{\pm 1}(c) &\equiv& \dim\, H_j(N^-\cup\{c\},N^- )^{\pm\Z_m}. \lb{2.13}
\eea
Clearly the integers $k_j(c)$ and $k_j^{\pm 1}(c)$ equal to
$0$ when $j<0$ or $j>\nu(c)$ and can take only values $0$ or $1$
when $j=0$ or $j=\nu(c)$.}

{\bf Lemma 2.4.} (cf. Lemma 2.4 of \cite{LoW1}) {\it Let $c$ be a prime
closed geodesic on a Finsler manifold $(M, F)$ satisfying (Iso).

(i) There holds
\bea 0\le k_j^{\pm 1}(c^m) \le k_j(c^m), \qquad
           \forall\,m\in\N,\;j\in\Z.  \lb{2.14}\eea

(ii) For any $m\in\N$, there holds
\bea k_0^{+1}(c^m) = k_0(c^m),\quad k_0^{-1}(c^m) = 0. \lb{2.15}\eea

(iii) In particular, if $c^m$ is non-degenerate, i.e.
$\nu(c^m)=0$, then
\bea k_0^{+1}(c^m) = k_0(c^m)=1,\quad k_0^{-1}(c^m) = 0. \lb{2.16}\eea}

Following the ideas of Gromoll-Meyer on the degenerate part of the
critical module in Theorem 3 of \cite{GrM2}, we have the following
result.

{\bf Proposition 2.5.} (cf. Theorem 3 of \cite{GrM2}, Section 7.1
of \cite{Rad2} and Theorem 3.11 of \cite{BaL1}) {\it Let $c$ be a
prime closed geodesic on a Finsler manifold $(M, F)$ satisfying
(Iso). Suppose for some integer $m=np\ge 2$ with $n$ and $p\in\N$
the nullitities satisfy
$$ \nu(c^m)=\nu(c^n).$$
Then the following holds for the degenerate part of the critical
module of $E$ with coefficient $\Q$.

(i) For any integer $j$, there hold
\bea
&& H_j(N_{c^m}^-\cup\{c^m\},N_{c^m}^-)
       = H_j(N_{c^n}^-\cup\{c^n\},N_{c^n}^-), \nn\\
&& k_j(c^m)=k_j(c^n). \lb{2.17}\eea

(ii) For any integer $j$, there hold
\bea
H_j(N_{c^m}^-\cup\{c^m\},N_{c^m}^-)^{\pm\Z_m}
&=& H_j(N_{c^n}^-\cup\{c^n\},N_{c^n}^-)^{\pm\Z_n}, \nn\\
{k}_j^{\pm 1}(c^m) &=& {k}^{\pm 1}_j(c^n).  \lb{2.18}\eea }

{\bf Proposition 2.6.} (cf. Satz 6.13 of \cite{Rad2}) {\it Let $c$
be a prime closed geodesic on a Finsler manifold $(M, F)$
satisfying (Iso). For any $m\in\N$, we have

(i) If $k_0(c^m)=1$, there holds $k_j^{\pm 1}(c^m)=0$ for
$1\le j\le \nu(c^m)$.

(ii) If $k_{\nu(c^m)}^{+1}(c^m)=1$ or $k_{\nu(c^m)}^{-1}(c^m)=1$,
there holds $k_j^{\pm 1}(c^m)=0$ for $0\le j\le \nu(c^m)-1$.

(iii) If $k_j^{+1}(c^m)\ge 1$ or $k_j^{-1}(c^m)\ge 1$ for some
$1\le j\le \nu(c^m)-1$, there holds $k_{\nu(c^m)}^{\pm 1}(c^m)=0=k_0(c^m)$.

(iv) In particular, if $\nu(c^m)\le 2$, then only one of the
$k_j(c^m)$'s can be non-zero.}

{\bf Proposition 2.7.}  {\it Let $c$ be a
prime closed geodesic on a Finsler manifold $(M, F)$ satisfying
(Iso). Suppose for some integer $m=np\ge 2$ with $n$ and $p\in\N$
the nullities satisfy
$$ \nu(c^m)\ge\nu(c^n). \lb{2.38}$$
Suppose further $k_{\nu(c^m)}(c^m)=1$, then we have  $k_{\nu(c^n)}(c^n)=1$.
In particular, we have
\bea k_j(c^n)=0,\quad 0\le j\le\nu(c^n)-1.\lb{2.19}\eea }

{\bf Proof.} The proof follows directly from Lemma 5 of \cite{GrM1}
as well as Section 7.1 of \cite{Rad2}, we now describe it briefly.
Let $D_{c^m}$ be a $\Z_m$-invariant
local hypersurface transverse to $S^1\cdot c^m$ in $c^m\in D_{c^m}$
as above and similarly for $D_{c^n}$. Then the $p$-fold iteration map
$\phi^p$ maps $D_{c^n}$ into $D_{c^m}$.
Since  $E|D_{c^m}$ is $\Z_p$-invariant, $\grad E(d^p)$ is targential to the
fixed point set $\phi^p(D_{c^n})$ for $d\in D_{c^n}$ and this yields
\bea \grad E(d^p)=\phi_\ast^p(\grad E(d)),\quad \forall d\in D_{c^n}.\lb{2.20}\eea
Now the proof of Lemma 1 in \cite{GrM1} yields that $\phi^p$ embeds
$N_{c^n}$  into $N_{c^m}$ as a submanifold, where $N_{c^n}$ and $N_{c^m}$
are local characteristic manifolds  at $c^n$ and $c^m$ respectively.
Note that
\be E(\phi^p(d))=p^2E(d), \quad \forall d\in N_{c^n}. \lb{2.21}\ee
By Corollary 8.4 of \cite{MaW1}, $k_{\nu(c^m)}(c^m)=1$
if and only if $c^m$ is a local maximum of $E$  in $N_{c^m}$.
This together with (\ref{2.21}) imply that  $c^n$ is a local
maximum of $E$  in $N_{c^n}$ too. Hence we use Corollary 8.4 \
of \cite{MaW1} again to obtain the proposition.\hfill\hb

\setcounter{equation}{0}
\section{The structure of $H_{\ast}(\ol{\Lm} S^n,\ol{\Lm}^0S^n;\Q)$ }

In this section, we briefly describe the relative homological
structure of the quotient space $\overline{\Lm}\equiv
\overline{\Lm} S^n$. Here we have $\ol{\Lm}^0=\ol{\Lambda}^0S^n
=\{{\rm constant\;point\;curves\;in\;}S^n\}\cong S^n$.

{\bf Theorem 3.1.} (H.-B. Rademacher, Theorem 2.4 and Remark 2.5 of \cite{Rad1})
{\it We have the Poincar\'e series

(i) When $n=2k+1$ is odd
\bea
P(\ol{\Lm}S^n,\ol{\Lm}^0S^n)(t)
&=&t^{n-1}\left(\frac{1}{1-t^2}+\frac{t^{n-1}}{1-t^{n-1}}\right)\nn\\
&=& t^{2k}\left(\frac{1}{1-t^2}+\frac{t^{2k}}{1- t^{2k}}\right).
\lb{3.1}\eea
Thus for $q\in\Z$ and $l\in\N_0$, we have
\bea {b}_q &=& {b}_q(\ol{\Lm}S^n,\ol{\Lm}^0 S^n)\nn\\
 &=&\rank H_q(\ol{\Lm} S^n,\ol{\Lm}^0 S^n )\nn\\
 &=& \;\;\left\{\matrix{
    2,&\quad {\it if}\quad q\in \{4k+2l,\quad l=0\;\mod\; k\},  \cr
    1,&\quad {\it if}\quad q\in \{2k\}\cup\{2k+2l,\quad l\neq 0\;\mod\; k\},  \cr
    0 &\quad {\it otherwise}. \cr}\right. \lb{3.2}\eea

(ii) When $n=2k$ is even
\bea
P(\ol{\Lm}S^n,\ol{\Lm}^0S^n)(t)
&=&t^{n-1}\left(\frac{1}{1-t^2}+\frac{t^{n(m+1)-2}}{1-t^{n(m+1)-2}}\right)
\frac{1-t^{nm}}{1-t^n}\nn\\
&=& t^{2k-1}\left(\frac{1}{1-t^2}+\frac{t^{4k-2}}{1- t^{4k-2}}\right),
\lb{3.3}\eea
where $m=1$ by Theorem 2.4 of \cite{Rad1}. Thus for $q\in\Z$ and $l\in\N_0$, we have}
\bea {b}_q &=& {b}_q(\ol{\Lm}S^n,\ol{\Lm}^0 S^n)\nn\\
 &=&\rank H_q(\ol{\Lm} S^n,\ol{\Lm}^0 S^n )\nn\\
 &=& \;\;\left\{\matrix{
    2,&\quad {\it if}\quad q\in \{6k-3+2l,\quad l=0\;\mod\; 2k-1\},  \cr
    1,&\quad {\it if}\quad q\in \{2k-1\}\cup\{2k-1+2l,\quad l\neq 0\;\mod\; 2k-1\},  \cr
    0 &\quad {\it otherwise}. \cr}\right. \lb{3.4}\eea

We have the following version of the Morse inequality.

{\bf Theorem 3.2.} (Theorem 6.1 of \cite{Rad2}) {\it Suppose that there exist
only finitely many prime closed geodesics $\{c_j\}_{1\le j\le p}$ on $(M, F)$,
and $0\le a<b\le \infty$ are regular values of the energy functional $E$.
Define for each $q\in\Z$,
\bea
{M}_q(\ol{\Lm}^b,\ol{\Lm}^a)
&=& \sum_{1\le j\le p,\;a<E(c^m_j)<b}\rank{\ol{C}}_q(E, c^m_j ) \nn\\
{b}_q(\ol{\Lm}^{b},\ol{\Lm}^{a})
&=& \rank H_q(\ol{\Lm}^{b},\ol{\Lm}^{a}). \nn\eea
Then there holds }
\bea
M_q(\ol{\Lm}^{b},\ol{\Lm}^{a}) &-& M_{q-1}(\ol{\Lm}^{b},\ol{\Lm}^{a})
    + \cdots +(-1)^{q}M_0(\ol{\Lm}^{b},\ol{\Lm}^{a}) \nn\\
&\ge& b_q(\ol{\Lm}^{b},\ol{\Lm}^{a}) - b_{q-1}(\ol{\Lm}^{b},\ol{\Lm}^{a})
   + \cdots + (-1)^{q}b_0(\ol{\Lm}^{b},\ol{\Lm}^{a}), \lb{3.5}\\
{M}_q(\ol{\Lm}^{b},\ol{\Lm}^{a}) &\ge& {b}_q(\ol{\Lm}^{b},\ol{\Lm}^{a}).\lb{3.6}
\eea

\setcounter{equation}{0}
\section{Classification of closed geodesics on $S^n$ and existence theorem}

Let $c$ be a closed geodesic on a Finsler n-sphere $S^n=(S^n,\,F)$.
Denote the linearized Poincar\'e map of $c$ by $P_c\in\Sp(2n-2)$.
Then $P_c$ is a symplectic matrix.
Note that the index iteration formulae in \cite{Lon1} of 2000 (cf. Chap. 8 of
\cite{Lon2}) work for Morse indices of iterated closed geodesics (cf.
\cite{LLo1}, Chap. 12 of \cite{Lon2}). Since every closed geodesic
on a  sphere must be orientable. Then by Theorem 1.1 of \cite{Liu1}
of C. Liu (cf. also \cite{Wil1}), the initial Morse index of a closed geodesic
$c$ on a $n$-dimensional Finsler  sphere coincides with the index of a
corresponding symplectic path introduced by C. Conley, E. Zehnder, and Y. Long
in 1984-1990 (cf. \cite{Lon2}).

As in \S 1.8 of \cite{Lon2}, define the homotopy component $\Omega^0 (P_c)$
of $P_c$ to be the path component of $\Omega (P_c)$, where
\bea \Omega (P_c) =\{N\in
Sp(2n-2) \mid &&\sigma (N)\cap U=\sigma (P_c)\cap U, \;and\nn\\
&&\nu_\lambda(N)=\nu_\lambda(P_c) \;\forall \lambda\in \sigma
(P_c)\cap U\}.\lb{4.1}\eea
The next theorem is due to Y. Long (cf. Theorem 8.3.1 and Corollary 8.3.2 of \cite{Lon2}).

{\bf Theorem 4.1.} {\it Let $\gamma\in\{\xi\in
C([0,\tau],Sp(2n))\mid \xi(0)=I\}$, Then there exists a path $f\in
C([0,1],\Omega^0(\gamma(\tau))$ such that $f(0)=\gamma(\tau)$ and
\bea f(1)=&&N_1(1,1)^{\diamond p_-} \diamond I_{2p_0}\diamond
N_1(1,-1)^{\diamond p_+}
\diamond N_1(-1,1)^{\diamond q_-} \diamond (-I_{2q_0})\diamond
N_1(-1,-1)^{\diamond q_+}\nn\\
&&\diamond R(\theta_1)\diamond\cdots\diamond R(\theta_r)
\diamond N_2(\omega_1, u_1)\diamond\cdots\diamond N_2(\omega_{r_*}, u_{r_*}) \nn\\
&&\diamond N_2(\lm_1, v_1)\diamond\cdots\diamond N_2(\lm_{r_0}, v_{r_0})
\diamond M_0 \lb{4.2}\eea
where $ N_2(\omega_j, u_j) $s are
non-trivial and   $ N_2(\lm_j, v_j)$s  are trivial basic normal
forms; $\sigma (M_0)\cap U=\emptyset$; $p_-$, $p_0$, $p_+$, $q_-$,
$q_0$, $q_+$, $r$, $r_*$ and $r_0$ are non-negative integers;
$\omega_j=e^{\sqrt{-1}\alpha_j}$, $
\lambda_j=e^{\sqrt{-1}\beta_j}$; $\theta_j$, $\alpha_j$, $\beta_j$
$\in (0, \pi)\cup (\pi, 2\pi)$; these integers and real numbers
are uniquely determined by $\gamma(\tau)$. Then using the
functions defined in (\ref{1.2}).
\bea i(\gamma, m)=&&m(i(\gamma,
1)+p_-+p_0-r)+2\sum_{j=1}^r\mathcal{E}\left(\frac{m\theta_j}{2\pi}\right)-r
-p_--p_0\nn\\&&-\frac{1+(-1)^m}{2}(q_0+q_+)+2\left(
\sum_{j=1}^{r_*}\varphi\left(\frac{m\alpha_j}{2\pi}\right)-r_*\right).
\lb{4.3}\eea
\bea \nu(\gamma, m)=&&\nu(\gamma,
1)+\frac{1+(-1)^m}{2}(q_-+2q_0+q_+)+2(r+r_*+r_0)\nn\\
&&-2\left(\sum_{j=1}^{r}\varphi\left(\frac{m\theta_j}{2\pi}\right)+
\sum_{j=1}^{r_*}\varphi\left(\frac{m\alpha_j}{2\pi}\right)
+\sum_{j=1}^{r_0}\varphi\left(\frac{m\beta_j}{2\pi}\right)\right)\lb{4.4}\eea
\bea \hat i(\gamma, 1)=i(\gamma, 1)+p_-+p_0-r+\sum_{j=1}^r
\frac{\theta_j}{\pi}.\lb{4.5}\eea
Where $N_1(1, \pm 1)=
\left(\matrix{ 1 &\pm 1\cr 0 & 1\cr}\right)$, $N_1(-1, \pm 1)=
\left(\matrix{ -1 &\pm 1\cr 0 & -1\cr}\right)$,
$R(\theta)=\left(\matrix{\cos\th &
                  -\sin\th\cr\sin\th & \cos\th\cr}\right)$,
$ N_2(\omega, b)=\left(\matrix{R(\th) & b
                  \cr 0 & R(\th)\cr}\right)$    with some
$\th\in (0,\pi)\cup (\pi,2\pi)$ and $b=
\left(\matrix{ b_1 &b_2\cr b_3 & b_4\cr}\right)\in\R^{2\times2}$,
such that $(b_2-b_3)\sin\theta>0$, if $ N_2(\omega, b)$ is
trivial; $(b_2-b_3)\sin\theta<0$, if $ N_2(\omega, b)$ is
non-trivial. We have $i(\gamma, 1)$ is odd if $f(1)=N_1(1, 1)$, $I_2$,
$N_1(-1, 1)$, $-I_2$, $N_1(-1, -1)$ and $R(\theta)$; $i(\gamma, 1)$ is
even if $f(1)=N_1(1, -1)$ and $ N_2(\omega, b)$; $i(\gamma, 1)$ can be any
integer if $\sigma (f(1)) \cap \U=\emptyset$.}

We will use the following theorem of N. Hingston in the $S^n$
case below. Note that the proof of N. Hingston's theorem does
not use the special properties of Riemannian metric, hence it
holds for Finsler metric as well.

{\bf Theorem 4.2.} (Follows from Proposition 1 of \cite{Hin1}, cf.
Lemma 3.4.12 of \cite{Kli3})
{\it Let $c$ be a closed geodesic of length $L$ on a Finsler
$n$-sphere $S^n=(S^n,F)$ such that as a critical orbit of the
energy functional $E$ on $\Lm S^n$, every orbit $S^1\cdot c^m$ of
its iteration $c^m$ is isolated. Suppose
\bea
&& i(c^m)+\nu(c^m) \le m(i(c)+\nu(c))-(n-1)(m-1),  \qquad
                \forall m\in\N, \lb{4.6}\\
&& k_{\nu(c)}(c)\neq0.  \lb{4.7} \eea
Then  $S^n$ has infinitely many prime closed geodesics. }

Note that in (\ref{4.7}), we have used the Shifting theorem in
\cite{GrM1}. Especially, (\ref{4.7}) means that $c$ is a local
maximum in the local characteristic manifold $N_c$ at $c$.

\setcounter{equation}{0}
\section{A mean index equality on $(S^n, F)$}

In this section, suppose that there are only finitely many prime closed
geodesics $\{c_j\}_{1\le j\le p}$ on $(S^n, F)$ with $\hat i(c_j)>0$
for $1\le j\le p$. We establish an equality of $c_j$s involving their
mean indices.

{\bf Definition 5.1.} {\it Let $M=(M, F)$ be a compact Finsler manifold
of dimension $n$, $c$ be a prime closed geodesic on $M$ satisfying (Iso).
For each $m\in\N$, the critical type numbers of $c^m$ is defined by the
following $2n-1$ tuple of integers via Definition 2.3
\bea K(c^m)
&\equiv&(k_0^\beta(c^m), k_1^\beta(c^m), \ldots, k_{2n-2}^\beta(c^m))\nn\\
&=&(k_0^{\beta (c^m)}(c^m),  k_1^{\beta(c^m)}(c^m), \ldots,
     k_{\nu (c^m)}^{\beta (c^m)}(c^m), 0, \ldots, 0), \lb{5.1}\eea
where $\beta=\bb(c^m)=(-1)^{i(c^m)-i(c)}$. We call a prime closed
geodesic $c$ homologically invisible if $K(c^m)=0$ for all $m\in\N$ and
homologically visible otherwise.}

As Lemma 2 of \cite{GrM2}, we have

{\bf Lemma 5.2.} {\it  Let $c$ be a prime closed geodesic on a
compact Finsler manifold $(M,F)$ satisfying (Iso).
Then there exists a minimal $T(c)\in \N$ such that
\bea
&& \nu(c^{p+T(c)})=\nu(c^p),\quad i(c^{p+T(c)})-i(c^p)\in 2\Z,
     \qquad\forall p\in \N,  \lb{5.2}\\
&& K(c^{p+T(c)})=K(c^p), \qquad\forall p\in \N. \lb{5.3}\eea
We call $T(c)$ the minimal period of critical modules of iterations of $c$. }

{\bf Proof.} Denote the linearized Poincar\'e
map of $c$ by $P_c:\R^{2(n-1)}\to\R^{2(n-1)}$. Then $P_c$ is a
symplectic matrix.
Denote by $\lambda_i=exp(\pm2\pi\frac{r_i}{s_i})$
the eigenvalues of $P_c$ possessing rotation angles which are rational
multiple of $\pi$ with $r_i$, $s_i\in\N$ and $(r_i,s_i)=1$ for
$1\le i\le q$. Let $T(c)$ be the least common multiple of
$s_1,\ldots, s_q$. Then the first equality in (\ref{5.2}) holds. If the
second equality in (\ref{5.2}) does not hold, replace $T(c)$ by $2T(c)$.
Then the later conclusion in (\ref{5.2}) follows from Theorem 9.3.4
of \cite{Lon2}.

In order to prove (\ref{5.3}), it suffices to show
\be K(c^{m+qT(c)})=K(c^m), \qquad \forall q\in \N,
         \; 1\le m\le T(c).  \lb{5.4}\ee

In fact, assume that (\ref{5.4}) is proved. Note that (\ref{5.3}) follows
from (\ref{5.4}) with $q=1$ directly when $p\le T(c)$. When $p>T(c)$, we
write $p=m+qT(c)$ for some $q\in\N$ and $1\le m\le T(c)$. Then by (\ref{5.4})
we obtain
$$ K(c^{p+T(c)})=K(c^{m+(q+1)T(c)})=K(c^m)=K(c^{m+qT(c)})=K(c^p), $$
i.e., (\ref{5.3}) holds.

To prove (\ref{5.4}), we fix an integer $m\in [1,T(c)]$. Let
$$ A=\{s_i\in\{s_1,\ldots, s_q\}\;|\;\;s_i\;{\rm is\;a\;factor\;of\;}m\},  $$
and let $m_1$ be the least common multiple of elements in $A$. Hence we have
$m=m_1m_2$ for some $m_2\in \N$ and $\nu(c^m)=\nu(c^{m_1})$.
Thus by Proposition 2.5 we have $k_l^{\beta(c^m)}(c^m)=k_l^{\beta(c^m)}(c^{m_1})$.
Since $m+pT(c)=m_1m_3$ for some $m_3\in\N$, we have by Proposition 2.5 that
$k_l^{\beta(c^{m+pT(c)})}(c^{m+pT(c)})=k_l^{\beta(c^{m+pT(c)})}(c^{m_1})$.
By (\ref{5.2}), we obtain $\beta(c^{m+pT(c)})=\beta(c^m)$, and then
(\ref{5.4}) is proved. This completes the proof. \hfill\hb

{\bf Definition 5.3.} {\it The Euler characteristic $\chi(c^m)$
of $c^m$ is defined by
\bea \chi(c^m)
&\equiv& \chi\left((\Lm(c^m)\cup S^1\cdot c^m)/S^1, \Lm(c^m)/S^1\right), \nn\\
&\equiv& \sum_{q=0}^{\infty}(-1)^q\dim \overline{C}_q( E,c^m).
\lb{5.5}\eea
Here $\chi(A, B)$ denotes the usual Euler characteristic of the space pair $(A, B)$.

The average Euler characteristic $\hat\chi(c)$ of $c$ is defined by }
\be \hat{\chi}(c)=\lim_{N\to\infty}\frac{1}{N}\sum_{1\le m\le N}\chi(c^m).
\lb{5.6}\ee

The following remark shows that $\hat\chi(c)$ is well-defined and is a
rational number.

{\bf Remark 5.4.} By (\ref{5.5}), we have \be \chi(c^m) =
\sum_{q=0}^{\infty}(-1)^q\dim \overline{C}_q( E,c^m) =
\sum_{l=0}^{2n-2}(-1)^{i(c^m)+l}k_l^{\beta(c^m)}(c^m).\lb{5.7}\ee
Here the  second equality follows from Proposition 2.1 and
Definition 5.1. Hence by (\ref{5.6}) and Lemma 5.2 we have \bea
\hat\chi(c) &=&\lim_{N\rightarrow\infty}\frac{1}{N}
  \sum_{1\le m\le N\atop 0\le l\le 2n-2}(-1)^{i(c^m)+l}k_l^{\beta(c^m)}(c^m)\nn\\
&=&\lim_{s\rightarrow\infty}\frac{1}{sT(c)}
  \sum_{1\le m\le T(c),\; 0\le l\le 2n-2\atop 0\le p< s}
  (-1)^{i(c^{pT(c)+m})+l}k_l^{\beta(c^{pT(c)+m})}(c^{pT(c)+m})\nn\\
&=&\frac{1}{T(c)}
  \sum_{1\le m\le T(c)\atop 0\le l\le 2n-2}
  (-1)^{i(c^{m})+l}k_l^{\beta(c^m)}(c^{m})
=\frac{1}{T(c)}\sum_{1\le m\le T(c)}\chi(c^m).\lb{5.8}\eea
Therefore $\hat\chi(y)$ is well defined and is a rational number.

Let $(X,Y)$ be a space pair such that the Betti numbers
$b_i=b_i(X,Y)=\dim H_i(X,Y;\Q)$ are finite for all $i\in \Z$. As usual
the {\it Poincar\'e series} of $(X,Y)$ is defined by the formal power
series $P(X, Y)=\sum_{i=0}^{\infty}b_it^i$. Suppose there exist only
finitely many  prime closed geodesics $\{c_j\}_{1\le j\le p}$
and  satisfy $\hat i(c_j)>0$ for $1\le j\le p$ on $(S^n, F)$.
The {\it Morse series} $M(t)$ of the energy functional $E$ of the space
pair $(\Lambda S^n/S^1, \Lambda^0 S^n/S^1)$ is defined as usual by
$$ M(t)=\sum_{q\ge 0,\; m\ge 1 \atop 1\le j\le p}\dim\ol{C}_q(E, c^m_j)t^q. $$
Then Theorem 3.2 yields a formal power series $Q(t)=\sum_{i=0}^\infty q_it^i$
with nonnegative integer coefficients $q_i$ such that
\be  M(t)=P(\Lm S^n/S^1,\Lm^0 S^n/S^1)(t)+(1+t)Q(t). \lb{5.9}\ee
For a formal power series $R(t)=\sum_{i=0}^\infty r_it^i$, we denote
by $R^n(t)=\sum_{i=0}^n r_i t^i$ for $n\in\N$ the corresponding
truncated polynomials. Using this notation, (\ref{5.9}) becomes
\be (-1)^Iq_I=M^I(-1)-P^I(-1), \qquad \forall I\in\N.    \lb{5.10}\ee
By Satz 7.8 of \cite{Rad2} we have
\bea &&\lim_{I\to\infty}\frac{1}{I}P^I(\Lm S^n/S^1,\Lm^0 S^n/S^1)(-1)\nn\\
=&&B(n,m)\nn\\
=&&\left\{\matrix{
     \frac{-m(m+1)n}{2n(m+1)-4}, &\quad n\quad even, \cr
     \frac{n+1}{2(n-1)} &\quad n\quad odd, \cr}\right.
     \lb{5.11}\eea
where $m=1$ by Corollary 2.6 of \cite{Rad1}.

The following consequence of an important result by H.-B. Rademacher (Theorem 7.9
in \cite{Rad2}) is needed in Sections 6 below. Here we have derived
precise dependence of coefficients on prime closed geodesics in the mean
index equality in Theorem 7.9 of \cite{Rad2}. This precise dependence is
also crucial for our proofs in Sections 6  below.

{\bf Theorem 5.5.} {\it Suppose that there exist only finitely many prime
closed geodesics $\{c_j\}_{1\le j\le p}$ with $\hat i(c_j)>0$
for $1\le j\le p$ on $(S^n,F)$.
Then the following identity holds
\be  \sum_{1\le j\le p}\frac{\hat\chi(c_j)}{\hat{i}(c_j)}=B(n,1).\lb{5.12}\ee
}

{\bf Proof.} Because $\dim \ol{C}_q(E,c_j^m)$ can be non-zero only for
$q=i(c_j^m)+l$ with $0\le l\le 2n-2$ by Proposition 2.1,
the formal Poincar\'e series $M(t)$ becomes
\be M(t)=\sum_{1\le j\le p,\; 0\le l\le 2n-2 \atop m\ge 1}
              k_l^{\beta}(c^m_j)t^{i(c_j^m)+l}
  = \sum_{1\le j\le p,\; 0\le l\le 2n-2 \atop 1 \le m\le T_j,\; s\ge 0}
               k_l^{\beta}(c_j^m)t^{i(c_j^{sT_j+m})+l}, \lb{5.13}\ee
where we denote by $T_j=T(c_j)$ for $1\le j\le p$.
The last equality follows from Lemma 5.2. Write
$M(t)=\sum_{h=0}^{\infty}w_ht^h$. Then we have
\be w_h\ = \sum_{1\le j\le p,\; 0\le l\le 2n-2 \atop 1 \le m\le T_j}
              k_l^{\beta}(c_j^m)\,^\#\{s\in\N_0\,|\,i(c_j^{sT_j+m})+l=h\}.
     \lb{5.14}\ee

{\bf Claim 1.} {\it $\{w_h\}_{h\ge 0}$ is bounded}.

In fact, we have
\bea
^\#\{s\in\N_0 &|& i(c_j^{sT_j+m})+l=h \}\nn\\
&=&\;^\#\{s\in\N_0 \;| \; i(c_j^{sT_j+m})+l=h,\;
                          |i(c_j^{sT_j+m})-(sT_j+m)\hat{i}(c_j)|\le n-1\} \nn\\
&\le &\;^\#\{s\in\N_0 \;| \;|h-l-(sT_j+m)\hat{i}(c_j)|\le n-1\}  \nn\\
&=&\;^\#\left\{s\in\N_0 \; \left|\;\frac{}{}\right.
      \;\frac{h-l-(n-1)-m\hat{i}(c_j)}{T_j\hat{i}(c_j)}\le s
       \le \frac{h-l+(n-1)-m\hat{i}(c_j)}{T_j\hat{i}(c_j)}\right\}  \nn\\
&\le&\; \frac{2(n-1)}{T_j\hat{i}(c_j)}+1,  \nn\eea
where the first equality follows from the fact $|i(c^m)-m\hat{i}(c)|\le n-1$
(cf. Theorem 1.4 of \cite{Rad1}). Hence Claim 1 holds.

We estimate next $M^I(-1)$. By (\ref{5.8}) we obtain
\bea M^I(-1) &=& \sum_{h=0}^I w_h(-1)^h   \nn\\
&=& \sum_{1\le j\le p,\; 0\le l\le 2n-2 \atop 1 \le m\le T_j}
            (-1)^{i(c_j^m)+l}k_l^{\beta}(c_j^m)
              \,^\#\{s\in\N_0 \,|\, i(c_j^{sT_j+m})+l\le I\}. \lb{5.15}\eea
Here the latter  equality holds by Lemma 5.2.

{\bf Claim 2.} {\it There is a real constant $C>0$ independent of $I$,
but depend on $c_j$ for $1\le j\le p$  such that
\be \left|M^I(-1)-\sum_{1\le j\le p,\; 0\le l\le 2n-2 \atop 1 \le m\le T_j}
 (-1)^{i(c_j^m)+l}k_l^{\beta}(c_j^m)\frac{I}{T_j\hat{i}(c_j)}\right|
               \le C,   \lb{5.16}\ee
where the sum in the left hand side of (\ref{5.16}) equals to
$I\sum_{1\le j\le p }\frac{\hat\chi(c_j)}{\hat i(c_j)}$ by (\ref{5.8}). }

In fact, we have the estimates
\bea
^\#\{s\in\N_0 &|& i(c_j^{sT_j+m})+l\le I\}   \nn\\
&=&\;^\#\{s\in\N_0 \;| \; i(c_j^{sT_j+m})+l\le I,\;
             |i(c_j^{sT_j+m})-(sT_j+m)\hat{i}(c_j)|\le n-1\}  \nn\\
&\le&\;^\#\{s\in\N_0 \;| \;0\le (sT_j+m)\hat{i}(c_j)\le I-l+(n-1)\}  \nn\\
&=&\;^\#\left\{s\in\N_0 \; \left |\;\frac{}{}\right.
    \;0\le s\le \frac{I-l+(n-1)-m\hat{i}(c_j)}{T_j\hat{i}(c_j)}\right\}  \nn\\
&\le&\; \frac{I-l+(n-1)}{T_j\hat{i}(c_j)}+1.  \nn\eea
On the other hand, we have
\bea
^\#\{s\in\N_0 &|& i(c_j^{sT_j+m})+l\le I\}  \nn\\
&=&\;^\#\{s\in\N_0 \;| \; i(c_j^{sT_j+m})+l\le I,\;
               |i(c_j^{sT_j+m})-(sT_j+m)\hat{i}(c_j)|\le n-1\}  \nn\\
&\ge&\;^\#\{s\in\N_0\;|\;i(c_j^{sT_j+m})\le(sT_j+m)\hat{i}(c_j)+n-1\le I-l\} \nn\\
&\ge&\;^\#\left\{s\in\N_0 \; \left |\;\frac{}{}\right.
   \;0\le s\le \frac{I-l-(n-1)-m\hat{i}(c_j)}{T_j\hat{i}(c_j)}\right\}  \nn\\
&\ge&\;\frac{I-l-(n-1)}{T_j\hat{i}(c_j)}-2,  \nn\eea
where $m\le T_j$ is used.
Combining these two estimates together with (\ref{5.15}), we obtain (\ref{5.16}).

By Claim 1, the sequence $\{w_h\}_{h\ge 0}$ is bounded. Hence by (\ref{5.9}),
the coefficient sequence $\{q_h\}_{h\ge 0}$ of $Q(t)$ is
bounded. Dividing both sides of (\ref{5.10}) by $I$, and letting $I$ tend
to infinity, together with Claim 2 and (\ref{5.11}) we obtain
$$ \lim_{I\to\infty}\frac{1}{I}M^I(-1)
          =\lim_{I\to\infty}\frac{1}{I}P^I(-1)=B(n,1). $$
Hence (\ref{5.12}) holds by (\ref{5.16}).  \hfill\hb

\setcounter{equation}{0}
\section{ Stability of closed geodesics on $(S^n,\, F)$ }

In this section, we give the proofs of Theorems 1.1, 1.2 and 1.4 by using
the mean index identity in Theorem 5.5, Morse inequality and the index iteration theory developed by
Y. Long and his coworkers.

In the following for the notation introduced in Section 3 we use
specially $M_j=M_j(\ol{\Lm} S^n,\ol{\Lm}^0 S^n)$ and
$b_j=b_j(\ol{\Lm} S^n,\ol{\Lm}^0 S^n)$ for $j=0,1,2,\ldots$.

{\bf Proof of Theorem 1.1.} First note that if the flag
curvature $K$ of $(S^n, F)$ satisfies
$\left(\frac{\lambda}{\lambda+1}\right)^2<K\le 1$,
then every nonconstant closed geodesic must satisfy
\bea i(c)\ge n-1. \lb{6.1}\eea
This follows from Theorem 3 and Lemma 3 of \cite{Rad3}.
Now it follows from Theorem 2.2 of \cite{LoZ1}
(Theorem 10.2.3 of \cite{Lon2}) that
\bea i(c^{m+1})-i(c^m)-\nu(c^m)\ge i(c)-\frac{e(P_c)}{2}\ge 0,\quad\forall m\in\N.\lb{6.2}\eea
Here the last inequality holds by (\ref{6.1}) and the fact that $e(P_c)\le 2(n-1)$.

Next we prove the theorem by showing that:
If the number of  prime closed geodesics is finite, then there must
exist at least one elliptic closed geodesic whose linearized
Poincar\'e map has at least one eigenvalue which is an irrational
multiple of $\pi$.

In the rest of this paper, we will assume the following

{\bf (F) There are only finitely many prime closed geodesics
$\{c_j\}_{1\le j\le p}$ on $(S^n,\,F)$. }

The proof contains two steps. In the first step, we prove the existence of
at least one elliptic closed geodesic. In the second step, we use Theorem 4.2
and the ideas of step 1 to complete the proof of the theorem.

{\bf Step 1.} {\it Claim: Under the assumption (F),
there must exist at least one elliptic closed geodesic. }

We prove it by contradiction, i.e., suppose that $e(P_{c_j})<2n-2$ for
$1\le j\le p$, where $P_{c_j}$ denotes the linearized Poincar\'e map
of $c_j$. Since $e(P_{c_j})$ is always even, we have
\bea e(P_{c_j})\le 2n-4,\quad 1\le j\le p.\lb{6.3}\eea

Note that by (\ref{6.1}) and (\ref{4.5}), we have $\hat i(c_j)>0$ for
$1\le j\le p$. Actually, we have $\hat i(c_j)> n-1$ for $1\le j\le p$
under the pinching assumption by Lemma 2 of \cite{Rad4}. Hence we can use the common index jump theorem (Theorem 4.3 of
\cite{LoZ1}, Theorem 11.2.1 of \cite{Lon2}) to obtain some
$(N, m_1,\ldots,m_p)\in\N^{p+1}$ such that
\bea
i(c_j^{2m_j}) &\ge& 2N-\frac{e(P_{c_j})}{2}, \lb{6.4}\\
i(c_j^{2m_j})+\nu(c_j^{2m_j}) &\le& 2N+\frac{e(P_{c_j})}{2}, \lb{6.5}\\
i(c_j^{2m_j+1}) &=& 2N+i(c_j). \lb{6.6}\eea
Moreover $\frac{m_j\theta}{\pi}\in\Z$, whenever $e^{\sqrt{-1}\theta}\in\sigma(P_{c_j})$
and $\frac{\theta}{\pi}\in\Q$. More precisely, by Theorem 4.1 of
\cite{LoZ1} (in (11.1.10) in Theorem 11.1.1 of \cite{Lon2}, with $D_j=\hat i(c_j)$,
we have
\bea m_j=\left(\left[\frac{N}{M\hat i(c_j)}\right]+\xi_j\right)M,\quad 1\le j\le p,\lb{6.7}\eea
where $\xi_j=0$ or $1$ for $1\le j\le p$ and $\frac{M\theta}{\pi}\in\Z$,
whenever $e^{\sqrt{-1}\theta}\in\sigma(P_{c_j})$ and $\frac{\theta}{\pi}\in\Q$
for some $1\le j\le p$. Furthermore, by (11.1.20) in Theorem 11.1.1 of \cite{Lon2},
for any $\epsilon>0$, we can choose $N$ and $\{\xi_j\}_{1\le j\le p}$ such that
\bea \left|\frac{N}{M\hat i(c_j)}-\left[\frac{N}{M\hat i(c_j)}\right]-\xi_j\right|<\epsilon
<\frac{1}{1+\sum_{1\le j\le p}4M|\hat\chi(c_j)|},\quad 1\le j\le p.\lb{6.8}\eea
Now by (\ref{6.1})-(\ref{6.6}), we have
\bea
i(c_j^{m})+\nu(c_j^{m}) &\le& i(c_j^{2m_j}),\quad\forall m<2m_j, \lb{6.9}\\
i(c_j^{2m_j})+\nu(c_j^{2m_j}) &\le& 2N+n-2, \lb{6.10}\\
i(c_j^{m}) &\ge& 2N+n-1,\quad\forall m>2m_j. \lb{6.11}\eea
By Theorem 5.5 and (\ref{5.11}), we have
\be  \sum_{1\le j\le p}\frac{\hat\chi(c_j)}{\hat{i}(c_j)}=B(n,1)\in\Q.\lb{6.12}\ee
Note by the proof of Theorem 4.1 of \cite{LoZ1} (Theorem 11.1.1 of \cite{Lon2}),
we can require that $N\in\N$  further satisfies (cf. (11.1.22) in \cite{Lon2})
\bea 2NB(n,1)\in\Z.\lb{6.13}\eea
Multiplying both sides of (\ref{6.12}) by $2N$ yields
\be  \sum_{1\le j\le p}\frac{2N\hat\chi(c_j)}{\hat{i}(c_j)}=2NB(n,1).\lb{6.14}\ee

{\bf Claim 1.} {\it We have}
\be  \sum_{1\le j\le p}2m_j\hat\chi(c_j)=2NB(n,1).\lb{6.15}\ee
In fact, by (\ref{6.12}), we have
\bea &&2NB(n,1)\nn\\
=&&\sum_{1\le j\le p}\frac{2N\hat\chi(c_j)}{\hat{i}(c_j)}\nn\\
=&&\sum_{1\le j\le p}2\hat\chi(c_j)\left(\left[\frac{N}{M\hat i(c_j)}\right]+\xi_j\right)M
+\sum_{1\le j\le p}2\hat\chi(c_j)\left(\frac{N}{M\hat i(c_j)}-\left[\frac{N}{M\hat i(c_j)}\right]-\xi_j\right)M\nn\\
=&& \sum_{1\le j\le p}2m_j\hat\chi(c_j)+\sum_{1\le j\le p}2M\hat\chi(c_j)\epsilon_j.\lb{6.16}
\eea
By Lemma 5.2 and our choice of $M$, we have
\bea \frac{2m_j}{T(c_j)}\in\N,\quad 1\le j\le p.\lb{6.17}\eea
Hence (\ref{5.8}) implies that
\bea 2m_j\hat\chi(c_j)\in\Z,\quad1\le j\le p.\lb{6.18}\eea
Now Claim 1 follows by (\ref{6.8}), (\ref{6.16}), (\ref{6.13}) and  (\ref{6.18}).

{\bf Claim 2.} {\it  We have
\be \sum_{1\le j\le p}2m_j\hat\chi(c_j)=M_0-M_1+M_2-\cdots+(-1)^{2N+n-2}M_{2N+n-2}.\lb{6.19}\ee}
In fact, by definition, the right hand side of (\ref{6.19}) is
\bea RHS=\sum_{q\le 2N+n-2\atop m\ge 1,\; 1\le j\le p}(-1)^q\dim\ol{C}_q(E, c^m_j).\lb{6.20}
\eea
By (\ref{6.9})-(\ref{6.11}) and Proposition 2.1, we have
\bea RHS&=&\sum_{1\le j\le p,\;1\le m\le 2m_j\atop q\le 2N+n-2}(-1)^q\dim\ol{C}_q(E, c^m_j),\lb{6.21}\\
&=&\sum_{1\le j\le p,\;1\le m\le 2m_j}\chi(c_j^m),\lb{6.22}
\eea
where the second equality follows from (\ref{5.7}), (\ref{6.9})-(\ref{6.10})
and Proposition 2.1.

By Lemma 5.2, (\ref{5.7})-(\ref{5.8}) and (\ref{6.17}), we have
\bea  \sum_{1\le m\le 2m_j}\chi(c_j^m)
&=&\sum_{0\le s< 2m_j/T(c_j)\atop1\le m\le T(c_j)}\chi(c_j^{sT(c_j)+m})\nn\\
&=&\frac{2m_j}{T(c_j)}\sum_{1\le m\le T(c_j)}\chi(c_j^m)\nn\\
&=& 2m_j\hat\chi(c_j),\lb{6.23}\eea
This proves Claim 2.

In order to complete Step 1, we have to consider the following two
cases according to the parity of $n$.

{\bf Case 1. } $n=2k+1$ is odd.

In this case, we have by (\ref{5.11})
\be B(n,1)=\frac{n+1}{2(n-1)}=\frac{k+1}{2k}.\lb{6.24}\ee
By the proof of Theorem 4.1 of \cite{LoZ1} (Theorem 11.1.1 of \cite{Lon2}),
we may  further assume  $N=mk$ for some $m\in\N$.

Thus by (\ref{6.15}), (\ref{6.19}) and (\ref{6.24}), we have
\be M_0-M_1+M_2-\cdots+(-1)^{2N+n-2}M_{2N+n-2}=m(k+1).\lb{6.25}
\ee
On the other hand, we have by (\ref{3.2})
\bea &&b_0-b_1+b_2-\cdots+(-1)^{2N+n-2}b_{2N+n-2}\nn\\
=&&b_{2k}+(b_{2k+2}+\cdots+b_{4k}+\cdots+b_{2mk+2}+\cdots+b_{2mk+2k})-b_{2mk+2k}\nn\\
=&& 1+m(k-1+2)-2\nn\\
=&& m(k+1)-1.\lb{6.26}
\eea
In fact, we cut off the sequence $\{b_{2k+2},\ldots,b_{2mk+2k}\}$ into
$m$ pieces, each of them contains  $k$ terms. Moreover, each piece contain
$1$ for $k-1$ times and  $2$ for one time. Thus (\ref{6.26}) holds.

Now by Theorem 3.2 and (\ref{6.26}), we have
\bea -m(k+1)&=&M_{2N+n-2}-M_{2N+n-3}+\cdots+M_1-M_0\nn\\
&\ge&b_{2N+n-2}-b_{2N+n-3}+\cdots+b_1-b_0\nn\\
&=&-(m(k+1)-1).\lb{6.27}
\eea
This contradiction yields Step 1 for $n$ being odd.

{\bf Case 2. } $n=2k$ is even.

In this case, we have by (\ref{5.11})
\be B(n,1)=\frac{-2n}{4n-4}=\frac{-k}{2k-1}.\lb{6.28}\ee
As in Case 1, we may assume $N=m(2k-1)$ for some $m\in\N$.

Thus by (\ref{6.15}), (\ref{6.19}) and (\ref{6.28}), we have
\be M_0-M_1+M_2-\cdots+(-1)^{2N+n-2}M_{2N+n-2}=-2mk.\lb{6.29}
\ee
On the other hand, we have by (\ref{3.4})
\bea &&b_0-b_1+b_2-\cdots+(-1)^{2N+n-2}b_{2N+n-2}\nn\\
=&&-b_{2k-1}-(b_{2k+1}+\cdots+b_{6k-3}+\cdots+b_{(m-1)(4k-2)+2k+1}+\cdots+b_{m(4k-2)+2k-1})\nn\\
&&+b_{m(4k-2)+2k-1}\nn\\
=&& -1-m(2k-2+2)+2\nn\\
=&& -2mk+1.\lb{6.30}
\eea
In fact, we cut off the sequence $\{b_{2k+1},\ldots,b_{m(4k-2)+2k-1}\}$ into
$m$ pieces, each of them contains  $2k-1$ terms. Moreover, each piece contain
$1$ for $2k-2$ times and  $2$ for one time. Thus (\ref{6.30}) holds.

Now by (\ref{6.29})-(\ref{6.30}) and Theorem 3.2, we have
\bea -2mk&=&M_{2N+n-2}-M_{2N+n-3}+\cdots+M_1-M_0\nn\\
&\ge&b_{2N+n-2}-b_{2N+n-3}+\cdots+b_1-b_0\nn\\
&=&-2mk+1.\lb{6.31}
\eea
This contradiction yields Step 1 for $n$ being even.

{\bf Step 2.} Under the assumption (F), there must exist  one
elliptic closed geodesic whose linearized Poincar\'e map has at
least one eigenvalue which is of the form $\exp(\pi i \mu)$ with an
irrational $\mu$.

In fact, we shall prove a more stronger result. Denote by
$\{P_{c_j}\}_{1\le j\le p}$ the linearized Poincar\'e maps of them.
Suppose $\{M_{c_j}\}_{1\le j\le p}$ are the basic normal  form decompositions
of $\{P_{c_j}\}_{1\le j\le p}$ in
$\{\Omega^0(P_{c_j})\}_{1\le j\le p}$ as in Theorem 4.1. Then we have the following

{\bf Claim 3.} {\it  There must exist $d\in\{c_j\}_{1\le j\le p}$ such that the
following hold:

(i) $e(P_{d})=2n-2$, i.e., $d$ is elliptic.

(ii) $M_d$ does not contain $N_1(1,1)$s, $N_1(-1,-1)$s  and nontrivial $N_2(\omega, b)$s.

(iii) Any trivial $N_2(\omega, b)$ contained in  $M_d$ must satisfies
$\frac{\theta}{\pi}\in\Q$, where $\omega=e^{\sqrt{-1}}\theta$.

(iv) If  $M_d$ contains $R(\theta)$ for some $\frac{\theta}{\pi}\notin\Q$,
then $M_d$ does not contain $R(2\pi-\theta)$.

(v) $k_{\nu(d^{T(d)})}^{\beta(d^{T(d)})}(d^{T(d)})\neq 0$.
Hence $d^{T(d)}$ is a local maximum of the energy functional
in the local characteristic manifold at $d^{T(d)}$.

(vi) $M_d$ must contain a term $R(\theta)$ with $\frac{\theta}{\pi}\notin\Q$.
}

In fact, we first show that there must exist $d$ satisfying (i)-(iv).
Suppose none of the closed geodesics in $\{c_j\}_{1\le j\le p}$ satisfies
all of (i)-(iv). Then as in Step 1, we obtain some
$(N, m_1,\ldots,m_p)\in\N^{p+1}$ such that (\ref{6.4})-(\ref{6.6}) hold.

By Step 1, we have found an elliptic closed geodesic $c$ for which (\ref{6.3})
does not hold anymore. Our following argument is to find other conditions to
replace (\ref{6.3}), then use the proof of Step 1.

From (\ref{6.2}) and (\ref{6.4})-(\ref{6.6}), we have
\bea i(c_j^{2m_j})+\nu(c_j^{2m_j})&\le& i(c_j^{2m_j+1})-i(c_j)+\frac{e(P_{c_j})}{2},\lb{6.32}\\
&=& 2N+\frac{e(P_{c_j})}{2},\quad 1\le j\le p.\lb{6.33}\eea
Now if $c_j$ does not satisfy all of (ii)-(iv), then by the
proof of Theorem 2.2 in \cite{LoZ1}, The strict inequality in
(\ref{6.32}) must hold. In fact, (ii) follows from Cases 1, 3, 7 and
(iii) follows from  Case 8 of Theorem 2.2 in \cite{LoZ1} respectively.
Hence we only need to check the case (iv).
As in the proof of Theorem 2.2 in \cite{LoZ1}, it suffices to show that
\bea\nu(\gamma, m)-\frac{e(M)}{2}<i(\gamma, m+1)-i(\gamma,m)-i(\gamma, 1),\quad \forall m\in\N,\lb{6.34}\eea
when $M=\gamma(\tau)=R(\theta)\diamond R(2\pi-\theta)$ with $\frac{\theta}{\pi}\notin\Q$.
By Case 6 in P.336 of \cite{LoZ1} and the symplectic additivity of indices and nullities,
we have
\bea\nu(\gamma, m)-\frac{e(M)}{2}=2-2\varphi\left(\frac{m\theta}{2\pi}\right
)-2\varphi\left(\frac{m(2\pi-\theta)}{2\pi}\right)=-2.\lb{6.35}\eea
On the other hand, we have
\bea &&i(\gamma, m+1)-i(\gamma,m)-i(\gamma, 1)\nn\\
=&&2\left(E\left(\frac{(m+1)\theta}{2\pi}\right)+E\left(\frac{(m+1)(2\pi-\theta)}{2\pi}\right)\right)\nn\\
&&-2\left(E\left(\frac{m\theta}{2\pi}\right)+E\left(\frac{m(2\pi-\theta)}{2\pi}\right)\right)-2\nn\\
=&&0.\lb{6.36}\eea
In the last equality, we have used the fact that $E(a)+E(-a)=1$ whenever
$a\in(0,\,+\infty)\setminus\Z$.
Hence (\ref{6.34}) is true.
Now  (\ref{6.9})-(\ref{6.11}) still hold.
Hence the same proof as in Step 1 yields a contradiction.

We then show that there must exist $d$ satisfying (i)-(v).
Suppose none of the closed geodesics in $\{c_j\}_{1\le j\le p}$ satisfies
all of (i)-(v). Then it is easy to see that (\ref{6.15}),
(\ref{6.19})-(\ref{6.23}) still hold.
In fact, we only need to check (\ref{6.22}). We have
$$\overline{C}_{q}(E,\,c_j^m)=0,\qquad\forall m\le 2m_j,\;\; q\ge 2N+n-1,\;\;1\le j\le p,$$
which follows easily from Propositiom 2.1, (\ref{6.2}) and (11.2.4) in Theorem 11.2.1 of \cite{Lon2}
$$i(c_j^{2m_j-1})+\nu(c_j^{2m_j-1})=2N-(i(c_j)+2S^+_{P_{c_j}}(1)-\nu(c_j))\le 2N,$$
where we have used (\ref{6.1}) and the fact that $2S^+_{P_{c_j}}(1)-\nu(c_j)\ge -(n-1)$,
which follows from (15.4.21) in p.340 of \cite{Lon2}. This yields (\ref{6.22}).

Thus the same proof as in Step 1 yields a contradiction.

At last we prove that there must exist $d$ satisfying (i)-(vi).
Suppose none of the closed geodesics in $\{c_j\}_{1\le j\le p}$
satisfies all of (i)-(vi). Then by the above argument, we assume
$c_j$ for some $1\le j\le p$ satisfies all (i)-(v) but not (vi).
We consider $g=c_j^{T(c_j)}$. Then by (i)-(iii) and the assumption,
$P_g$ can be connected in
$\Omega^0(P_g)$ to $I_{2p_0^\prime}\diamond N_1(1,\,-1)^{\diamond p_+^\prime}$
with $p_0^\prime+p_+^\prime=n-1$ as in Theorem 4.1.
In fact, by (i)-(iii) in Claim 3 and the assumption, the basic normal
form decomposition (\ref{4.1}) in Theorem 4.1 becomes
\bea M_{c_j}=&&I_{2p_0}\diamond N_1(1,\,-1)^{\diamond p_+}
\diamond N_1(-1,\,1)^{\diamond q_-}\diamond(-I_{2q_0})\nn\\
&&\diamond R(\theta_1)\diamond\cdots\diamond R(\theta_r)
\diamond N_2(\lambda_1,\,v_1)\diamond\cdots\diamond N_2(\lambda_{r_0},\,v_{r_0})
\nn\eea
together with $\frac{\theta_j}{\pi}\in\Q$ for
$1\le j\le r$, $\frac{\beta_j}{\pi}\in\Q$  for $1\le j\le r_0$
and $p_0+p_++q_-+q_0+r+2r_0=n-1$.
By Lemma 5.2, we have $\frac{T(c_j)\theta_j}{2\pi}\in\Z$,
$\frac{T(c_j)\beta_j}{2\pi}\in\Z$  and $2|T(c_j)$ whenever
$-1\in\sigma(M_{c_j})$. Hence $R(\theta_j)^{T(c_j)}=I_2$,
$N_2(\lambda_j,\,v_j)^{T(c_j)}$ can be connected within
$\Omega^0(N_2(\lambda_j,\,v_j)^{T(c_j)})$ to $N_1(1,\,-1)^{\diamond 2}$
and $(-I_2)^{T(c_j)}=I_2$, $N_1(-1, \,1)^{T(c_j)}$
can be connected within $\Omega^0(N_1(-1, \,1)^{T(c_j)})$ to $N_1(1,\,-1)$
whenever $-1\in\sigma(M_{c_j})$. Thus $p_0^\prime=p_0+q_0+r$ and
$p_+^\prime=p_++q_-+2r_0$ and then $P_g$ behaves as claimed.

Now by Theorem 4.1, we have
\bea i(g^m)=m(i(g)+p_0^\prime)-p_0^\prime,\quad \nu(g^m)\equiv2p_0^\prime+p_+^\prime.\qquad\forall m\in\N.\lb{6.37}
\eea
Hence
\bea i(g^m)+\nu(g^m)=m(i(g)+p_0^\prime)+p_0^\prime+p_+^\prime.\qquad\forall m\in\N.\lb{6.38}
\eea
On the other hand
\bea &&m(i(g)+\nu(g))-(n-1)(m-1)\nn\\
=&&m(i(g)+2p_0^\prime+p_+^\prime)-(p_0^\prime+p_+^\prime)(m-1)\nn\\
=&&m(i(g)+p_0^\prime)+p_0^\prime+p_+^\prime.\qquad\forall m\in\N.\lb{6.39}
\eea
By Lemma 2.4 and (v), we have
\be k_{\nu(g)}(g)=k_{\nu(c_j^{T(c_j)})}(c_j^{T(c_j)})
\ge k_{\nu(c_j^{T(c_j)})}^{\beta(c_j^{T(c_j)})}(c_j^{T(c_j)})\neq 0.\lb{6.40}
\ee
Hence we can use Theorem 4.2 to obtain infinitely many prime closed geodesics,
which contradicts to the assumption (F). This complete the proof of Step 2.\hfill\hb

{\bf Proof of Theorem 1.2.} Suppose $d\in\{c_j\}_{1\le j\le p}$ is a
closed geodesic satisfying (i)-(vi) of Claim 3  above.
Then we have the following
\be i(d^m)+\nu(d^m)-i(d)-\nu(d)\in2\Z, \quad\forall m\in\N.\lb{6.41}\ee
One can prove this by verifying each basic norm form in the decomposition
of $M_d$, then using the symplectic additivity of indices and nullities
to obtain (\ref{6.41}). Here we omit the details.

By (v), we have $k_{\nu(d^{T(d)})}^{\beta(d^{T(d)})}(d^{T(d)})\neq 0$.
Hence by Lemma 2.4 and Proposition 2.7, we have
\bea k_{\nu(d^{T(d)})}(d^{T(d)})&\ge& k_{\nu(d^{T(d)})}^{\beta(d^{T(d)})}(d^{T(d)})=1,\lb{6.42}\\
0\le k_{\nu(d^m)}^{\beta(d^m)}(d^m)&\le& k_{\nu(d^m)}(d^m)=1,\qquad \forall m\in\N,\lb{6.43}\\
0\le k_l^{\beta(d^m)}(d^m)&\le& k_l(d^m)=0,\qquad \forall 0\le l\le \nu(d^m)-1,\; m\in\N.\lb{6.44}
\eea
Note that in order to get (\ref{6.43}) and  (\ref{6.44}),
we have used the same argument as in the last paragraph in the proof of Lemma 5.2.

By (\ref{5.8}), we have
\bea \hat\chi(d)
&=&\frac{1}{T(d)}\sum_{1\le m\le T(d)\atop 0\le l\le 2n-2}
  (-1)^{i(d^{m})+l}k_l^{\beta(d^m)}(d^{m})\lb{6.45}\\
&=&\frac{1}{T(d)}\sum_{1\le m\le T(d)}
  (-1)^{i(d^{m})+\nu(d^m)}k_{\nu(d^m)}^{\beta(d^m)}(d^{m})\lb{6.46}\\
&=&\frac{(-1)^{i(d)+\nu(d)}}{T(d)}\sum_{1\le m\le T(d)}
  k_{\nu(d^m)}^{\beta(d^m)}(d^{m})\lb{6.47}\\
&\neq& 0.\lb{6.48}\eea
Here to get (\ref{6.46}), we have used (\ref{6.44}).
To get (\ref{6.47}), we have used (\ref{6.41}).
To get (\ref{6.48}), we have used (\ref{6.42}) and (\ref{6.43}).
Note that $B(n, 1)\in\Q$. Hence the theorem follows
easily from (vi) of Claim 3, Theorem 4.1 and Theorem 5.5.\hfill\hb

{\bf Proof of Theorem 1.4.} Suppose $d\in\{c_j\}_{1\le j\le p}$ is a
closed geodesic satisfying (i)-(vi) of Claim 3  above.
Suppose $\omega^{\pm 1}_1,\ldots,\omega^{\pm 1}_r$ with $\omega_i=e^{\sqrt{-1}\theta_i}$
are those eigenvalues of $P_d$ that satisfy $\frac{\theta_i}{\pi}\notin\Q$ for $1\le i\le r$.
Then we have

Claim 4. {\it There exists $M\in \Sp(2n-2)$ such that
\be MP_dM^{-1}=R(\widehat{\theta_1})\diamond\cdots\diamond R(\widehat{\theta_r})\diamond M_0,\lb{6.49}\ee
with $\widehat{\theta_i}=\theta_i$ or $2\pi-\theta_i$ for $1\le i\le r$ and
$\sigma(M_0)\subset \U\cap\{e^{\sqrt{-1}\theta}\,|\,\frac{\theta}{\pi}\in\Q\}$.}

In fact, by Theorem 1.6.11 of \cite{Lon2}, we have $M_1\in \Sp(2n-2)$ such that
\be M_1P_dM_1^{-1}=S_1\diamond\cdots\diamond S_{m_1}\diamond S_0,\lb{6.50}\ee
with $S_0\in\Sp(2k_0)$ with $k_0\ge 0$ and $\omega_1\notin\sigma(S_0)$,
$k_i\ge 1$ and $S_i\in\Sp(2k_i)$ is of  the normal form $N_{k_i}(\lambda_i, b_i)$
with $\lambda_i=\omega_1$ or $\omega_1^{-1}$ defined in Section 1.6 of
\cite{Lon2}.

Then by Case 3 and 4 in Section 1.8 of \cite{Lon2},
if $k_i\ge 3$ for some $1\le i\le m_1$, then $S_i$ can be
connected within $\Omega^0(S_i)$ to $\widetilde{S_i}$ with
$e(\widetilde{S_i})<e(S_i)$. This contradicts to (i) of Claim 3.
Hence  $S_i=N_{l_i}(\lambda_i, b_i)$ for $l_i=1$ or $2$ and $1\le i\le m_1$.
We next prove that $l_i=1$ for $1\le i\le m_1$. Suppose $l_i=2$ for some $i$.
Then by (ii) and (iii) of Claim 3, $S_i$ is not a basic normal form,
hence $S_i\in \mathcal{M}_{\omega_1}^2(4)$. By Case 4 of \cite{Lon2},
$S_i$ can be connected within $\Omega^0(S_i)$ to
$R(\omega_1)\diamond R(2\pi-\omega_1)$. This contradicts to (iv) of Claim 3.
Now (\ref{6.50}) becomes
\be M_1P_dM_1^{-1}=R(\widehat{\theta_1})^{\diamond m_1}\diamond S_0,\lb{6.51}\ee
with $\widehat{\theta_1}=\theta_1$ or $2\pi-\theta_1$.
We continue the above argument for at most $r$ times and then obtain (\ref{6.49}).
This proves Claim 4.

At last, let $g=d^{T(d)}$. Then it is easy to see that $g$ is of
elliptic-parabolic type. In fact, we have
$$ MP_{d^{T(d)}}M^{-1}=(MP_{d}M^{-1})^{T(d)}=R(T(d)\widehat{\theta_1})\diamond\cdots\diamond R(T(d)\widehat{\theta_r})\diamond M_0^{T(d)}$$
while by the proof Lemma 5.2, all the eigenvalues of $M_0^{T(d)}$ equal to $1$.
Hence Theorem 1.4 holds by the definition of elliptic-parabolic type.\hfill\hb

\setcounter{equation}{0}
\section{ Existence of three closed geodesics on $(S^3,\, F)$ }

In this section, we give a proof of Theorem 1.5 based on the results
in the previous sections.

Firstly by Lemma 2 of \cite{Rad4}, we have the following

{\bf Lemma 7.1.} {\it Suppose the Finsler  $3$-sphere $(S^3,\,F)$
with reversibility $\lambda$ and flag curvature $K$
satisfies $\left(\frac{\lambda}{\lambda+1}\right)^2< K\le 1$. Then
every closed geodesic $c$ on $(S^3,\,F)$ has mean index $\hat
i(c)>2$. }

{\bf Remark 7.2.} Note that by Theorem 4.1, the index iteration
formulae of $N_1(1, \,1)$ and $I_2$ coincide and both can be viewed
as a rotation matrix $R(\theta)$ with $\theta=2\pi$.
Similarly $N_1(-1, \,-1)$ and $-I_2$  can be viewed
as a rotation matrix $R(\theta)$ with $\theta=\pi$.
Hence in the following, we only consider the case $R(\theta)$
with $\theta\in (0,\,2\pi]$, the same proof still works for
$N_1(1, \,1)$ and $N_1(-1, \,-1)$.

{\bf Lemma 7.3.} {\it Under the assumption of Lemma 7.1,
suppose $P_c$ can be connected within $\Omega^0(P_c)$ to $M_c$
as in Theorem 4.1 with $M_c=R(\theta_1)\diamond R(\theta_2)$ with
$\theta_j\in (0,\,2\pi]$ for $j=1, 2$. Then we have the following}
\bea i(c^{m+q})\ge i(c^m)+2,\qquad \forall m\in\N,\; q\ge 2.\lb{7.1}
\eea
{\bf Proof.} By Theorem 4.1, we have
\bea i(c^ m)=m(i(c)-2)+
2\mathcal{E}\left(\frac{m\theta_1}{2\pi}\right)
+2\mathcal{E}\left(\frac{m\theta_2}{2\pi}\right)-2,\quad\forall m\in\N,\lb{7.2}
\eea
with $\theta_j\in (0,\,2\pi]$ for $j=1, 2$ and $i(c)\in 2\N$.

Now if $i(c)\ge 4$, the lemma is obvious by (\ref{7.2}).

Next consider the case $i(c)=2$. We have
\bea i(c^ m)=
2\mathcal{E}\left(\frac{m\theta_1}{2\pi}\right)
+2\mathcal{E}\left(\frac{m\theta_2}{2\pi}\right)-2,\quad\forall m\in\N,\lb{7.3}
\eea
By Lemma 7.1, we have
$$ \hat i(c)=\frac{\theta_1}{\pi}+\frac{\theta_2}{\pi}>2.$$
Hence with out loss of generality, we may assume $\frac{\theta_1}{\pi}>1$.
Thus we have
\bea i(c^{m+p})&=&
2\mathcal{E}\left(\frac{(m+p)\theta_1}{2\pi}\right)
+2\mathcal{E}\left(\frac{(m+p)\theta_2}{2\pi}\right)-2\nn\\
&\ge&2\mathcal{E}\left(\frac{m\theta_1}{2\pi}+\frac{\theta_1}{\pi}\right)
+2\mathcal{E}\left(\frac{m\theta_2}{2\pi}\right)-2\nn\\
&\ge&2\mathcal{E}\left(\frac{m\theta_1}{2\pi}\right)
+2\mathcal{E}\left(\frac{m\theta_2}{2\pi}\right),\nn
\eea
for $p\ge 2$. Note that in the second inequality above, we have
used the fact that $\mathcal{E}(a+b)\ge \mathcal{E}(a)+1$
for any $a\in\R$ and $b>1$. Hence the lemma holds.\hfill\hb

{\bf Lemma 7.4.} {\it Under the assumption of Lemma 7.1,
suppose $P_c$ can be connected to $M_c$ as in Theorem 4.1 with
$M_c=R(\theta_1)\diamond R(\theta_2)$ with $\frac{\theta_j}{\pi}\in (0,\,2]\cap\Q$
for $j=1, 2$. Let $m\in\N$ satisfy $\nu(c^m)=4$. Then we have the following}
\bea i(c^{m-q})+\nu(c^{m-q})\le i(c^m)-2,\qquad \forall  q\ge 2.\lb{7.4}
\eea
{\bf Proof.} As in Lemma 7.3, we have (\ref{7.2}) with
\be\frac{\theta_j}{2\pi}=\frac{r_j}{s_j},\quad r_j,\,s_j\in\N,\;\;(r_j,\, s_j)=1,\;\; j=1,2.\lb{7.5}\ee
By (\ref{6.2}), it suffices to prove the case $q=2$.

Note that $\nu(c^m)=4$ implies that $s_j|m$, i.e., $s_j$ is a factor of $m$ for
$j=1,\,2$. We may assume $s_1\le s_2$ without loss of generality. Hence we have
\be \nu(c^{m-1})=\nu(c),\qquad\nu(c^{m-2})=\nu(c^2).\lb{7.6}\ee
In fact, $\nu(c^{m-1})=2k$ for $k\in\{0,\,1,\,2\}$ if and only if  $s_j|m-1$ for $1\le j\le k$,
and this is equivalent to $s_j|1$ for $1\le j\le k$, and this implies $\nu(c)=2k$.
Similarly, $\nu(c^{m-2})=2k$ for $k\in\{0,\,1,\,2\}$ if and only if  $s_j|m-2$ for $1\le j\le k$,
and this is equivalent to $s_j|2$ for $1\le j\le k$, and this implies $\nu(c^2)=2k$.

Now if $\nu(c)\ge2$, then by (\ref{6.2}) and (\ref{7.6}), we have
$$ i(c^{m-2})+\nu(c^{m-2})\le i(c^{m-1})\le i(c^m)-\nu(c^{m-1})\le i(c^m)-2.$$
Hence the lemma holds.

Next consider the case $\nu(c)=0$. If $\nu(c^2)=0$,
then by Lemma 7.3 and (\ref{7.6}), we have
$$ i(c^{m-2})+\nu(c^{m-2})=i(c^{m-2})\le i(c^m)-2.$$
Hence it remains to consider the case $\nu(c^2)\ge 2$ and $\nu(c)=0$.
This implies $s_1=2$ and then $r_1=1$ by (\ref{7.5}).
Now we have
\bea i(c^ m)&=&m(i(c)-2)+
2\mathcal{E}\left(\frac{m}{2}\right)
+2\mathcal{E}\left(\frac{m\theta_2}{2\pi}\right)-2\nn\\
&=&m(i(c)-2)+
2\mathcal{E}\left(\frac{m-2}{2}\right)
+2\mathcal{E}\left(\frac{(m-2)\theta_2}{2\pi}+\frac{\theta_2}{\pi}\right)\nn\\
&=&i(c^{m-2})+2+2(i(c)-2)+
2\left(\mathcal{E}\left(\frac{(m-2)\theta_2}{2\pi}+\frac{\theta_2}{\pi}\right)
-\mathcal{E}\left(\frac{(m-2)\theta_2}{2\pi}\right)\right).\lb{7.7}
\eea
Now if $i(c)\ge 4$, then (\ref{7.7}) implies that
$i(c^m)\ge i(c^{m-2})+6$.
This together with $\nu(c^{m-2})\le 4$ prove the lemma.

If $i(c)= 2$, then by Lemma 7.1, we have
$\hat i(c)=1+\frac{\theta_2}{\pi}>2$. This yields $\frac{\theta_2}{\pi}>1$,
and then $s_2\ge 3$. Hence  $\nu(c^{m-2})\le 2$ by (\ref{7.6}).
Now the last term in (\ref{7.7}) is not less than $2$. Hence
$$i(c^{m-2})+\nu(c^{m-2})\le i(c^m)-4+2=i(c^m)-2.$$
This proves the whole lemma.\hfill\hb

{\bf Lemma 7.5.} {\it Let $M_j$  and $b_j$ be the integers defined
at the beginning of Section 6. If $M_k=b_k$ for some $k\in\N_0$, then we have
\bea M_k-M_{k-1}+\cdots+(-1)^kM_0&=& b_k-b_{k-1}+\cdots+(-1)^kb_0,\lb{7.8}\\
M_{k-1}-M_{k-2}+\cdots+(-1)^{k-1}M_0&=& b_{k-1}-b_{k-2}+\cdots+(-1)^{k-1}b_0,\lb{7.9}
\eea}
{ \bf Proof.} By Theorem 3.2, we have
\bea M_k-M_{k-1}+\cdots+(-1)^kM_0&\ge& b_k-b_{k-1}+\cdots+(-1)^kb_0,\nn\\
M_{k-1}-M_{k-2}+\cdots+(-1)^{k-1}M_0&\ge& b_{k-1}-b_{k-2}+\cdots+(-1)^{k-1}b_0,\nn
\eea
These together with $M_k=b_k$ yields (\ref{7.9}) and then (\ref{7.8}).\hfill\hb

{\bf Lemma 7.6.} {\it Under the assumption of Lemma 7.1,
there are at least two prime closed geodesics on $(S^3,\, F)$.
If there are exactly two prime closed geodesics, then at least one of
them has Poincar\'e map $P_c=R(\theta_1)\diamond R(\theta_2)$ with
$\frac{\theta_j}{\pi}\in (0,\,2]\setminus\Q$ for $j=1, 2$ in an appropriate
coordinates. }

{\bf Proof.} By \cite{Fet1}, there exists at least one closed geodesic
on $(S^3,\, F)$. Now we prove $p\ge 2$, where $p$ is the integer in
the assumption (F). As in the proof of Theorem 1.1 in Section 6,
by Theorem 11.2.1 of \cite{Lon2}, we obtain some
$(N, m_1,\ldots,m_p)\in\N^{p+1}$ such that
\bea
i(c_j^{2m_j}) &\ge& 2N-\frac{e(P_{c_j})}{2}, \lb{7.10}\\
i(c_j^{2m_j})+\nu(c_j^{2m_j}) &\le& 2N+\frac{e(P_{c_j})}{2}, \lb{7.11}\\
i(c_j^{2m_j-m})+\nu(c_j^{2m_j-m}) &\le& 2N-(i(c_j)+2S^+_{P_{c_j}}(1)-\nu(c_j)),\quad \forall m\in\N. \lb{7.12}\\
i(c_j^{2m_j+m}) &\ge& 2N+i(c_j),\quad\forall m\in\N. \lb{7.13}\\
\nu(c_j^{2m_j+1})&=&\nu(c_j),\lb{7.14}
\eea
By the proof of Theorem 1.1, we may assume $c_1$ satisfies
\be M_{c_1}=R(\theta_1)\diamond Q_1,\quad i(c_1^{2m_1})+\nu(c_1^{2m_1})=2N+2,
\quad k_{\nu(c_1^{2m_1})}^{\beta(c_1^{2m_1})}(c_1^{2m_1})\neq 0,\lb{7.15}\ee
for some $\frac{\theta_1}{\pi}\in\R\setminus\Q$ and $Q_1\in Sp(2)$. Here we use
notations as in the proof of Theorem 1.1. In fact, by the proof of Theorem 1.1,
we must have $\overline{C}_{2N+2}(E,\, c_j^{2m_j})\neq 0$  and
$M_{c_j}$ contain a term $R(\theta)$ with $\frac{\theta}{\pi}\notin\Q$
for some $1\le j\le p$. Thus (\ref{7.15}) follows from Proposition 2.1.

By  p.340 of \cite{Lon2}, we have
\bea
&& 2S^+_{P_{c_1}}(1)-\nu(c_1) \nn\\
&&\qquad = 2S^+_{R(\theta_1)}(1)-\nu_1(R(\theta_1))
   +2S^+_{Q_1}(1)-\nu_1(Q_1)\nn\\
&&\qquad = 0 + 2S^+_{Q_1}(1) - \nu_1(Q_1)\nn\\
&&\qquad \ge -1.\lb{7.16}\eea
Now by (\ref{6.1}) we have
\bea i(c_1^{2m_1-m})+\nu(c_1^{2m_1-m}) &\le& 2N-1,\quad \forall m\in\N. \lb{7.17}\\
i(c_1^{2m_1+m}) &\ge& 2N+2,\quad\forall m\in\N. \lb{7.18}
\eea
Now (\ref{7.15}), (\ref{7.17})-(\ref{7.18}) together with Propositions 2.1 and 2.6
yield
\be \ol{C}_{2N}(E, c_1^m)=0,\quad\forall m\in\N.\lb{7.19}\ee
By Theorems 3.1 and 3.2 we have
\be M_{2N}\ge b_{2N}=2.\lb{7.20}\ee
Hence $p\ge 2$ holds. This proves the first part of the lemma.

Next we assume $p=2$. Then we have
\be
\ol{C}_{2N}(E, c_2^{q})\neq 0,\lb{7.21}\ee
for some $q\in\N$ by (\ref{7.20}).

In order to prove the second part of the lemma, we need the following

{\bf Claim.} The only possible value of $q$ that satisfies(\ref{7.21}) is  $2m_2$.

In fact, again by  p.340 of \cite{Lon2}, we have
$$2S^+_{P_{c_2}}(1)-\nu(c_2)\ge -2,$$
where the equality holds if and only if $M_{c_2}=N_1(1,\,-1)^{\diamond 2}$.
Now by (\ref{6.1}), (\ref{7.12}) and (\ref{7.13}),  we have
\bea i(c_2^{2m_2-m})+\nu(c_2^{2m_2-m}) &\le& 2N,\quad \forall m\in\N. \lb{7.22}\\
i(c_2^{2m_2+m}) &\ge& 2N+2,\quad\forall m\in\N. \lb{7.23}
\eea
Here the equality in (\ref{7.22})  holds if and only if  $m=1$,
$M_{c_2}=N_1(1,\,-1)^{\diamond 2}$ and $i(c_2)=2$ hold simultaneously.

If the claim is not true, then by (\ref{7.21}) and Proposition 2.1,
the equality for $m=1$ in (\ref{7.22}) must hold. Thus  by Theorem 4.1,
we have
$$ i(c_2^m)=2m,\quad \nu(c_2^m)=2,\qquad\forall m\in\N.$$
Hence by Proposition 2.1 again, we have
$$\ol{C}_{2N}(E, c_2^{2m_2-1})=\Q,\quad
\ol{C}_{2N}(E, c_2^{m})=0,\;\;\forall m\neq 2m_2-1.$$
This together with (\ref{7.19}) yield
$$ 1=M_{2N}\ge b_{2N}=2.$$
This contradiction proves the claim.

Now we prove the second part of the lemma by contradiction,
i.e., suppose $Q_1\neq R(\theta_2)$ with some
$\frac{\theta_2}{\pi}\in (0,\,2]\setminus\Q$ in the decomposition (\ref{7.15}).
Then by Theorem 1.2, we must have
\be M_{c_2}=R(\theta^\prime)\diamond Q^\prime,\lb{7.24}\ee
for some $\frac{\theta^\prime}{\pi}\in\R\setminus\Q$ and $Q^\prime\in Sp(2)$.

Now (\ref{7.19})-(\ref{7.21}) imply that
$$ M_{2N}=\dim\ol{C}_{2N}(E, c_2^{2m_2})\ge 2.$$
Note that by (\ref{7.24}), we have $\nu(c_2^{2m_2})\le 2$. If
$\nu(c_2^{2m_2})=2$, then by Theorem 4.1, we must have $Q^\prime=R(\vartheta)$ for
some $\frac{\vartheta}{\pi}\in (0,\,2]\cap\Q$. Thus it follows from
(\ref{4.3}) that $i(c_2^m)\in2\Z$ for all $m\in\N$. Hence either $2N=i(c_2^{2m_2})$ or
$2N=i(c_2^{2m_2})+2$ holds. Now by Definition 2.3 for $\nu(c_2^{2m_2})\le 2$, we have
$$\ol{C}_{2N}(E, c_2^{2m_2})=\Q.$$
This contradiction proves the lemma.\hfill\hb

Now we can give the proof of Theorem 1.5.

{\bf Proof of Theorem 1.5.} We prove the theorem by contradiction,
i.e., by Lemma 7.6, we assume $p=2$ in the assumption (F).
By Lemma 7.6, we may assume $P_{c_1}=R(\theta_1)\diamond R(\theta_2)$ with
$\frac{\theta_j}{\pi}\in (0,\,2]\setminus\Q$ for $j=1, 2$
and (\ref{7.15}) holds. Assume the linearized Poincar\'e map $P_{c_2}$ of $c_2$
can be connected to $M_{c_2}$ in Theorem 4.1. Then due to
Remark 7.2 and Theorem 4.1, we have the following
cases according to $M_{c_2}$.

{\bf Case 1.} {\it  We have $M_{c_2}=R(\vartheta_1)\diamond R(\vartheta_2)$ with
$\frac{\vartheta_j}{\pi}\in (0,\,2]\cap\Q$ for $j=1, 2$.}

By  p.340 of \cite{Lon2}, we have
\be2S^+_{P_{c_j}}(1)-\nu(c_j)\ge 0,\qquad j=1,\,2.\lb{7.25}\ee
Now by (\ref{6.1}), (\ref{7.12})-(\ref{7.13}), we have
\bea i(c_j^{2m_j-m})+\nu(c_j^{2m_j-m}) &\le& 2N-2,\quad \forall m\in\N,\lb{7.26}\\
i(c_j^{2m_j+m}) &\ge& 2N+2,\quad\forall m\in\N, \lb{7.27}
\eea
for $j=1,\,2$. Note that by the definition of $2m_2$, we have
$2m_2\frac{\vartheta_j}{2\pi}\in\Z$ for $j=1,\,2$. Thus we have
\be \nu(c_2^{2m_2})=4,\qquad i(c_2^{2m_2})=2N-2,\lb{7.28}\ee
where the latter holds by (\ref{7.10})-(\ref{7.11}).
By (\ref{7.2}), we have for $m\ge 2$
\bea i(c_2^{2m_2+m})&\ge &i(c_2^{2m_2})+2(i(c_2)-2)+
2\mathcal{E}\left(\frac{\vartheta_1}{\pi}\right)
+2\mathcal{E}\left(\frac{\vartheta_2}{\pi}\right)\nn\\
&\ge& i(c_2^{2m_2})+6\nn\\
&=&2N+4,\qquad \forall m\ge 2.\lb{7.29}
\eea
Where in the first inequality, we have used
$2m_2\frac{\vartheta_j}{2\pi}\in\Z$ for $j=1,\,2$.
The second inequality follows by $\hat i(c_2)>2$ as before.
Note that by (\ref{7.21}) and (iii) of Proposition 2.6, we have
\be\overline{C}_{2N-2}(E,\; c_2^{2m_2})=0=
\overline{C}_{2N+2}(E,\; c_2^{2m_2}).\lb{7.30}\ee
Hence (\ref{7.26})-(\ref{7.30}) imply
\be \sum_{m\in\N}\dim \overline{C}_{2N+2}(E,\; c_2^{m})\le 1.\lb{7.31}\ee
with equality holds if and only if $i(c_2)=2$ and
$k_0^{+1}(c_2^{2m_2+1})=k_0^{+1}(c_2)=1$. Note that here we have used $T(c_2)|2m_2$.

By (\ref{7.15}) and Lemma 7.3, we have
\be 1\le\sum_{m\in\N}\dim \overline{C}_{2N+2}(E,\; c_1^{m})\le 2.\lb{7.32}\ee
where the second equality holds if and only if $i(c_1)=2$.

By Theorem 3.2 and (\ref{7.31})-(\ref{7.32}), we have
$$2=b_{2N+2}\le M_{2N+2}=
\sum_{m\in\N,\;1\le j\le 2}\dim \overline{C}_{2N+2}(E,\; c_j^{m})\le 3.$$
Now if $M_{2N+2}=3$, by (\ref{7.31})-(\ref{7.32}), we have
\bea&&\overline{C}_{q}(E,\; c_1^m)=0=\overline{C}_{q}(E,\; c_2^m),
\qquad \forall m\in\N,\;q\le 1,\nn\\
&&\overline{C}_{2}(E,\; c_1)=\Q=\overline{C}_{2}(E,\; c_2),\nn\\
&&\overline{C}_{2}(E,\; c_1^m)=0=\overline{C}_{2}(E,\; c_2^m),\quad\forall m\ge 2,\nn\\
&&\overline{C}_{3}(E,\; c_1^m)=0=\overline{C}_{3}(E,\; c_2^m),\quad\forall m\in\N.\nn
\eea
Here the latter two hold by Propositions 2.1 and 2.6
together with $i(c_j^m)\ge 4$ for $m\ge 2$ and $j=1,\,2$,
which follows from $\hat i(c_j)>2$ for $j=1,\,2$ as before.
Hence by Theorems 3.1 and 3.2, we have
$$-2=M_3-M_2+M_1-M_0\ge b_3-b_2+b_1-b_0=-1.$$
This contradiction proves that
\be M_{2N+2}=b_{2N+2}=2.\lb{7.33}\ee
By (\ref{7.28}), Lemma 7.4, Proposition 2.1 and (\ref{7.30}), we have
\be \sum_{m\in\N}\dim \overline{C}_{2N-2}(E,\; c_2^{m})\le 1.\lb{7.34}\ee
with equality holds if and only if $i(c_2)=2$ and
$k_{\nu(c_2^{2m_2-1})}^{+1}(c_2^{2m_2-1})=k_{\nu(c_2)}^{+1}(c_2)=1$. Note that here we have used $T(c_2)|2m_2$.
By Lemma 7.3, we have
\be \sum_{m\in\N}\dim \overline{C}_{2N-2}(E,\; c_1^{m})\le 2.\lb{7.35}\ee
By Theorem 3.2 and (\ref{7.34})-(\ref{7.35}), we have
\be2=b_{2N-2}\le M_{2N-2}=
\sum_{m\in\N,\;1\le j\le 2}\dim \overline{C}_{2N-2}(E,\; c_j^{m})\le 3.\lb{7.36}\ee
Next we have two subcases according to  the value of $M_{2N-2}$.

{\bf Subcase 1.1.}  $M_{2N-2}=3$.

In this subcase, the equality in (\ref{7.34}) and (\ref{7.35}) must hold.
In particular, we have
\bea \overline{C}_{2N-3}(E,\; c_j^{m})=0,\qquad \forall m\in\N,\;\;j=1,\,2.\lb{7.37}\eea
In fact, since $c_1^m$ is non-degenerate and $i(c_1^m)$ is even for all $m\in\N$
by (\ref{7.2}). Hence (\ref{7.37}) holds for $c_1^m$ by Proposition 2.1.
By (\ref{7.28}), Lemma 7.4 and Proposition 2.1, (\ref{7.37}) holds for $c_2^m$
with $m\le 2m_2-2$. By (\ref{7.27})-(\ref{7.28}) and Proposition 2.1,
(\ref{7.37}) holds for $c_2^m$ with $m\ge 2m_2$.
Since the equality in (\ref{7.34}) holds, we have
$$i(c_2^{2m_2-1})+\nu(c_2^{2m_2-1})=2N-2,\qquad k_{\nu(c_2^{2m_2-1})}^{+1}(c_2^{2m_2-1})=1.$$
Hence (\ref{7.37}) holds for $c_2^{2m_2-1}$ by Propositions 2.1 and 2.6.

Now we have the following diagram.
\bea
\begin{tabular}{l|cccccc}
&$2N-3$&$2N-2$& $2N-1$&$2N$&$2N+1$&$2N+2$\\
\hline
$M_{\ast,\; 1}$ &$0$&$2$&$0$&$0$&$0$&$\xi$\\
$M_{\ast,\; 2}$ &$0$&$1$&$k_1$&$k_2$&$k_3$&$\zeta$\\
\end{tabular}\lb{7.38}\eea
where $M_{\ast,\; j}=\sum_{m\in\N}\dim \overline{C}_{\ast}(E,\; c_j^{m})$
for $j=1,\,2$, $k_l=k_l^{+1}(c_2^{2m_2})$ for $1\le l\le 3$ and $\xi+\zeta=2$.
In fact, the first column follows from (\ref{7.37}).
The second column follows from (\ref{7.34})-(\ref{7.35}).
The first row follows from (\ref{7.15}) and (\ref{7.26})-(\ref{7.27}).
The second row follows from  (\ref{7.26})-(\ref{7.27}).
$\xi+\zeta=2$ follows from  (\ref{7.33}).

Since $M_{2N-3}=b_{2N-3}=0$ and $M_{2N+2}=b_{2N+2}=2$ hold by the
first and last column of (\ref{7.38}), we have by Lemma 7.5
\bea M_{2N+2}-M_{2N+1}+\cdots-M_1+M_0&=& b_{2N+2}-b_{2N+1}+\cdots-b_1+b_0,\lb{7.39}\\
M_{2N-3}-M_{2N-4}+\cdots+M_1-M_0&=& b_{2N-3}-b_{2N-4}+\cdots+b_1-b_0,\lb{7.40}
\eea
Subtracting (\ref{7.40}) from (\ref{7.39}) and using Theorem 3.1 together
with (\ref{7.38}), we have
\bea 2-k_3+k_2-k_1+3&=&M_{2N+2}-M_{2N+1}+M_{2N}-M_{2N-1}+M_{2N-2}\nn\\
&=&b_{2N+2}-b_{2N+1}+b_{2N}-b_{2N-1}+b_{2N-2}=6.\nn
\eea
This implies
\be k_2-k_1-k_3=1.\lb{7.41}\ee
Since $c_1^m$ is non-degenerate and $i(c_1^m)$ is even for all $m\in\N$,
we have $T(c_1)=1$ by Lemma 5.2. Hence we have $\hat\chi(c_1)=1$ by (\ref{5.8}).
Now by Lemma 7.1, we have
\be \frac{\hat \chi(c_1)}{\hat i(c_1)}<\frac{1}{2}.\lb{7.42}\ee
Since $T(c_2)|2m_2$, we have
$k_l^{+1}(c_2^{T(c_2)})=k_l^{+1}(c_2^{2m_2})$ for all $l\in\N$.
This together with $i(c_2^m)$ is even for all $m\in\N$ and (\ref{7.41}), we have
\be \chi(c_2^{T(c_2)})=k_2-k_1-k_3=1.\lb{7.43}\ee
On the other hand, we have
\be \chi(c_2^m)\le 1,\qquad \forall m<T(c_2).\lb{7.44}\ee
In fact, this follows from  $i(c_2^m)$ is even, $\nu(c_2^m)\le 2$
for all $m<T(c_2)$, Proposition 2.6 and (\ref{5.7}),.
Now by (\ref{5.8}), we have
$$\hat\chi(c_2)=\frac{1}{T(c_2)}\sum_{1\le m\le T(c_2)}\chi(c_2^m)
\le\frac{1}{T(c_2)}T(c_2)=1.$$
Now by Lemma 7.1, we have
\be \frac{\hat \chi(c_2)}{\hat i(c_2)}<\frac{1}{2}.\lb{7.45}\ee
By Theorem 5.5, (\ref{5.11}), (\ref{7.42}) and (\ref{7.45}), we have
\be 1=B(3,\,1)=\frac{\hat \chi(c_1)}{\hat i(c_1)}
+\frac{\hat \chi(c_2)}{\hat i(c_2)}<\frac{1}{2}+\frac{1}{2}=1.\lb{7.46}\ee
This contradiction proves the theorem in this subcase.

{\bf Subcase 1.2.} $M_{2N-2}=2$.

In this subcase, we have the following diagram.
\bea
\begin{tabular}{l|ccccc}
&$2N-2$& $2N-1$&$2N$&$2N+1$&$2N+2$\\
\hline
$M_{\ast,\; 1}$ &$\eta$&$0$&$0$&$0$&$\xi$\\
$M_{\ast,\; 2}$ &$\lambda$&$k_1$&$k_2$&$k_3$&$\zeta$\\
\end{tabular}\lb{7.47}\eea
With  $\eta+\lambda=2$, $\xi+\zeta=2$ and
$k_l=k_l^{+1}(c_2^{2m_2})$ for $1\le l\le 3$.
Here $\eta+\lambda=2$ follows from $M_{2N-2}=2$ and the other parts follow
just as in Subcase 1.

Since $M_{2N-2}=b_{2N-2}=2$ and $M_{2N+2}=b_{2N+2}=2$ hold by the
first and last column of (\ref{7.47}), we have by Lemma 7.5
\bea M_{2N+2}-M_{2N+1}+\cdots-M_1+M_0&=& b_{2N+2}-b_{2N+1}+\cdots-b_1+b_0,\lb{7.48}\\
M_{2N-2}-M_{2N-3}+\cdots+M_1-M_0&=& b_{2N-2}-b_{2N-3}+\cdots+b_1-b_0,\lb{7.49}
\eea
Subtracting (\ref{7.49}) from (\ref{7.48}) and using Theorem 3.1 together (\ref{7.47}), we have
\bea 2-k_3+k_2-k_1&=&M_{2N+2}-M_{2N+1}+M_{2N}-M_{2N-1}\nn\\
&=&b_{2N+2}-b_{2N+1}+b_{2N}-b_{2N-1}=4.\nn
\eea
This implies
\be k_2-k_1-k_3=2.\lb{7.50}\ee
As in (\ref{7.43}), we have
\be \chi(c_2^{T(c_2)})=k_2-k_1-k_3=2.\lb{7.51}\ee
Now by (\ref{5.8}) and (\ref{7.44}), we have
\be\hat\chi(c_2)=\frac{1}{T(c_2)}\sum_{1\le m\le T(c_2)}\chi(c_2^m)
\le\frac{T(c_2)+1}{T(c_2)}.\lb{7.52}\ee
We assume
\be\frac{\vartheta_j}{2\pi}=\frac{r_j}{\lambda s_j},\quad r_j,\,\lambda,\,s_j\in\N,\;\;(r_j,\, \lambda s_j)=1,\;(s_1,\,s_2)=1,\; j=1,2.\lb{7.53}\ee
By (\ref{7.2}) and Lemma 7.1, we have
\bea \hat i(c_2)=i(c_2)-2+\frac{\vartheta_1}{\pi}+\frac{\vartheta_2}{\pi}
\equiv2q+\frac{2r_1}{\lambda s_1}+\frac{2r_2}{\lambda s_2}>2,\lb{7.54}
\eea
where we denote by $2q=i(c_2)-2\in2\N_0$.

Note that $T(c_2)=\lambda s_1s_2$. Thus Multiplying both sides of (\ref{7.54})
by $\frac{T(c_2)}{2}$ yields
$$ \frac{T(c_2)\hat i(c_2)}{2}=q\lambda s_1s_2+r_1s_2+r_2s_1>\lambda s_1s_2.$$
Hence
\be  \frac{T(c_2)\hat i(c_2)}{2}=q\lambda s_1s_2+r_1s_2+r_2s_1\ge \lambda s_1s_2+1.\lb{7.55}\ee
Now by (\ref{7.52}) and (\ref{7.55}), we have
\be \frac{\hat\chi(c_2)}{\hat i(c_2)}\le
\frac{T(c_2)+1}{T(c_2)\hat i(c_2)}=\frac{\lambda s_1s_2+1}{2(\lambda s_1s_2+1)}
=\frac{1}{2}.\lb{7.56}\ee
This together with (\ref{7.42}) yields the  contradiction (\ref{7.46}).
Hence the theorem holds in this subcase.

{\bf Case 2.} {\it  We have $M_{c_2}=R(\vartheta_1)\diamond N_1(-1,\,1)$ with
$\frac{\vartheta_1}{\pi}\in (0,\,2]\cap\Q$.}

By Theorem 4.1, we have
\be i(c_2^m)=m(i(c_2)-1)+2\mathcal{E}\left(\frac{m\vartheta_1}{2\pi}\right)-1,\quad\forall m\in\N,\lb{7.57}
\ee
with $i(c_2)\in 2\N$. Note that
\be i(c_2^m)-i(c_2)\in\left\{\matrix{2\N_0 & {\rm if}\; m\in2\N-1, \cr
2\N_0+1 & {\rm if}\; m\in2\N, \cr}\right. \lb{7.58}
\ee
One can check that (\ref{7.25})-(\ref{7.27}) still hold in this case.

Note that by the definition of $2m_2$, we have
$2m_2\frac{\vartheta_1}{2\pi}\in\Z$. Thus we have
\be \nu(c_2^{2m_2})=3,\qquad i(c_2^{2m_2})=2N-1,\lb{7.59}\ee
where the latter holds by (\ref{7.10})-(\ref{7.11}) and $i(c_2^{2m_2})$
is odd. Note that by (\ref{7.21}) and (iii) of Proposition 2.6, we have
\be\overline{C}_{2N-1}(E,\; c_2^{2m_2})=0=
\overline{C}_{2N+2}(E,\; c_2^{2m_2}).\lb{7.60}\ee
By (\ref{7.57}) and Lemma 7.1, we have
\be i(c_2^{2m_2+1})\ge 2N+2,\qquad i(c_2^{2m_2+m})\ge 2N+5,\quad\forall m\ge 2.\lb{7.61}\ee
Here the latter holds as in (\ref{7.29}). Thus by the same argument as in Case 1,
we have (\ref{7.33}).

Now we have the following diagram.
\bea
\begin{tabular}{l|cccc}
& $2N-1$&$2N$&$2N+1$&$2N+2$\\
\hline
$M_{\ast,\; 1}$ &$0$&$0$&$0$&$\xi$\\
$M_{\ast,\; 2}$ &$0$&$k_1$&$k_2$&$\zeta$\\
\end{tabular}\lb{7.62}\eea
With  $\xi+\zeta=2$ and $k_l=k_l^{-1}(c_2^{2m_2})$ for $1\le l\le 2$.
The first column follows by (\ref{7.15}), (\ref{7.26})-(\ref{7.27})
and (\ref{7.60}). The other parts follows just as in (\ref{7.38}).

Hence we have $M_{2N+2}=b_{2N+2}=2$ and $M_{2N-1}=b_{2N-1}=0$.
Thus as in Subcase 1.1, we obtain
\be 2-k_2+k_1=M_{2N+2}-M_{2N+1}+M_{2N}=b_{2N+2}-b_{2N+1}+b_{2N}=4.\lb{7.63}
\ee
This implies
\be k_1-k_2=2.\lb{7.64}\ee
Note that by Lemma 5.2 and (\ref{7.58}), $T(c_2)$ is even.
Hence by (\ref{7.57}), $T(c_2)|2m_2$ and $i(c_2^{2m_2})$ is odd, we have
\be \chi(c_2^{T(c_2)})=k_1-k_2=2.\lb{7.65}\ee
We assume
\be\frac{\vartheta_1}{2\pi}=\frac{r_1}{s_1},\quad r_1,\,s_1\in\N,\;\;(r_1,\, s_1)=1.\lb{7.66}\ee
By (\ref{7.57}) and Lemma 7.1, we have
\bea \hat i(c_2)=i(c_2)-1+\frac{\vartheta_1}{\pi}
\equiv q+\frac{2r_1}{s_1}>2,\lb{7.67}
\eea
where we denote by $q=i(c_2)-1\in2\N_0+1$.

Note that $T(c_2)=\frac{2s_1}{(2,\,s_1)}$. Thus Multiplying both sides of (\ref{7.67})
by $\frac{T(c_2)}{2}$ yields
$$ \frac{T(c_2)\hat i(c_2)}{2}=\frac{T(c_2)}{2}\left(q+\frac{2r_1}{s_1}\right) >T(c_2).$$
Since both the second and the third terms are integers, we have
\be  \frac{T(c_2)\hat i(c_2)}{2}\ge T(c_2)+1.\lb{7.68}\ee
Note that (\ref{7.44}) still holds. In fact, by Theorem 4.1,
we have $\nu(c_2^m)\le 1$ for $m\neq 0\,\mod \,s_1$. Hence ({\ref{7.44})
hold for these $m$ by Proposition 2.6. Now if $T(c_2)=s_1$, ({\ref{7.44})
is true. If $T(c_2)=2s_1$, i.e., $s_1$ is odd, we only need to consider
$\chi(c_2^{s_1})$, then $\nu(c_2^{s_1})=2$ and $i(c_2^{s_1})$ is even
by (\ref{7.58}). Thus $\chi(c_2^{s_1})\le 1$ holds by Proposition 2.6
and (\ref{5.7}).

By (\ref{7.44}), (\ref{5.8}), (\ref{7.65}) and (\ref{7.68}), we have
\be \frac{\hat\chi(c_2)}{\hat i(c_2)}\le
\frac{T(c_2)+1}{T(c_2)\hat i(c_2)}\le\frac{1}{2}.\lb{7.69}\ee
This together with (\ref{7.42}) yields the  contradiction (\ref{7.46}).
Hence the theorem holds in this case.

{\bf Case 3.} {\it  We have $M_{c_2}=R(\vartheta_1)\diamond N_1(1,\,-1)$ with
$\frac{\vartheta_1}{\pi}\in (0,\,2]\cap\Q$.}

By Theorem 4.1, we have
\be i(c_2^m)=m(i(c_2)-1)+2\mathcal{E}\left(\frac{m\vartheta_1}{2\pi}\right)-1,\quad\forall m\in\N,\lb{7.70}
\ee
with $i(c_2)\in 2\Z+1$. Hence by  (\ref{6.1}), we have $i(c_2)\ge 3$.
By  p.340 of \cite{Lon2}, we have $2S^+_{P_{c_2}}(1)-\nu(c_2)\ge -1$.
Thus by (\ref{7.12})-(\ref{7.13}), we have
\bea i(c_2^{2m_2-m})+\nu(c_2^{2m_2-m}) &\le& 2N-2,\quad \forall m\in\N,\lb{7.71}\\
i(c_2^{2m_2+m}) &\ge& 2N+3,\quad\forall m\in\N, \lb{7.72}
\eea
Note that the same argument as in Case 2 implies that (\ref{7.59})-(\ref{7.60}) still hold.
By (\ref{7.32}), Theorem 3.2, (\ref{7.60}), (\ref{7.71})-(\ref{7.72}) and Proposition 2.1,
we have (\ref{7.33}) here. By the same argument as in Case 2, (\ref{7.62}) holds here
with $k_l=k_l^{+1}(c_2^{2m_2})$ for $1\le l\le 2$. Then (\ref{7.63})-(\ref{7.64}) hold.

Note that $T(c_2)|2m_2$ and $i(c_2^{2m_2})$ is odd by (\ref{7.70}), we have
\be \chi(c_2^{T(c_2)})=k_1-k_2=2.\lb{7.73}\ee
We assume
\be\frac{\vartheta_1}{2\pi}=\frac{r_1}{s_1},\quad r_1,\,s_1\in\N,\;\;(r_1,\, s_1)=1.\lb{7.74}\ee
By (\ref{7.70}) and Lemma 7.1, we have
\bea \hat i(c_2)=i(c_2)-1+\frac{\vartheta_1}{\pi}
\equiv q+\frac{2r_1}{s_1}>2,\lb{7.75}
\eea
where we denote by $q=i(c_2)-1\in2\N$.

Note that $T(c_2)=s_1$. Thus Multiplying both sides of (\ref{7.67})
by $\frac{T(c_2)}{2}$ yields
$$ \frac{T(c_2)\hat i(c_2)}{2}=\frac{T(c_2)}{2}\left(q+\frac{2r_1}{s_1}\right) >T(c_2).$$
Note that $q\in 2\N$, hence both the second and the third terms are integers, hence we have
\be  \frac{T(c_2)\hat i(c_2)}{2}\ge T(c_2)+1.\lb{7.76}\ee
Note that (\ref{7.44}) still holds since $\nu(c_2^m)\le 1$ for $m<T(c_2)$.
Thus (\ref{7.44}), (\ref{7.73}) and (\ref{7.76}) imply (\ref{7.69}) holds.
This together with (\ref{7.42}) yields the  contradiction (\ref{7.46}).
Hence the theorem holds in this case.

{\bf Case 4.} {\it  We have $M_{c_2}=R(\vartheta_1)\diamond P$ with
$\frac{\vartheta_1}{\pi}\in (0,\,2]\cap\Q$ and $P\in Sp(2)$ is hyperbolic.}

In this case, by Theorem 4.1, the index iteration formula coincide with
that of Case 2 or Case 3 when $i(c_2)\in 2\N$ or $i(c_2)\in 2\N+1$ respectively.
Thus the same proof of Case 2 and 3 imply our theorem is true in this case.

{\bf Case 5.} {\it $M_{c_2}=R(\vartheta_1)\diamond R(\vartheta_2)$
with $\frac{\vartheta_1}{\pi}\in (0,\,2]\cap\Q$ and
$\frac{\vartheta_2}{\pi}\in (0,\,2]\setminus\Q$.}

In this case, by Theorem 4.1, we have
\be  i(c_2^{2m_2})\in2\N,\qquad \nu(c_2^{2m_2})=2.\lb{7.77}\ee
Then by the proof of Lemma 7.6 together with Propositions 2.1 and 2.6, we have
\be 1\ge\dim\ol{C}_{2N}(E, c_2^{2m_2})=M_{2N}\ge b_{2N}= 2.\lb{7.78}\ee
This contradiction proves the theorem in this case.

{\bf Case 6.} {\it $M_{c_2}$ dose not contain $R(\vartheta)$ with
$\frac{\vartheta}{\pi}\in (0,\,2]\cap\Q$ and $N_1(1,\,1)$, $N_1(-1,\,-1)$ .}

In this case, by Theorem 4.1, we have
$$ \nu(c_2^{2m_2})\le 2.$$
Moreover, if the equality holds, $M_{c_2}$ must be one of the following
$$ N_1(1,\,-1)^{\diamond 2},\quad N_1(1,\,-1)\diamond N_1(-1,\,1),\quad
N_1(-1,\,1)^{\diamond 2},\quad N_2(\omega,\,b),$$
with $N_2(\omega,\,b)=\left(\matrix{ R(\vartheta) & b\cr 0 & R(\vartheta)\cr}\right)$
and $\frac{\vartheta}{\pi}\in (0,\,2]\cap\Q$. Hence by Theorem 4.1, we have
$$ i(c_2^{2m_2})\in2\N.$$
Then the same argument as in Lemma 7.6 shows the theorem holds in this case.

Combining all the above cases, we obtain Theorem 1.5.\hfill\hb

\medskip

\noindent {\bf Acknowledgements.} I would like to sincerely thank the referee for his/her
careful reading and valuable comments and suggestions
on this paper. I would like to sincerely thank my
advisor, Professor Yiming Long, for introducing me to the theory of
closed geodesics and for his valuable help and encouragement during
the writing of this paper.

\bibliographystyle{abbrv}

\bigskip

\end{document}